\documentclass[11pt,a4paper]{article}

\usepackage[a4paper,margin=28mm]{geometry}
\usepackage{amsmath,amssymb,amsthm,mathtools,booktabs,array,tikz}
\makeatletter
\newenvironment{inplacefigure}{\par\addvspace{12pt}\noindent\begin{minipage}{\linewidth}\def\@captype{figure}\centering}{\end{minipage}\par\addvspace{12pt}}
\makeatother
\usepackage[hidelinks]{hyperref}
\usetikzlibrary{arrows.meta,calc}

\newcommand{\R}{\mathbb R}
\newcommand{\Pa}{\mathcal P}
\newcommand{\rank}{\operatorname{rank}}
\newcommand{\im}{\operatorname{im}}
\newcommand{\row}{\operatorname{row}}
\newcommand{\loc}{\mathrm{loc}}
\newcommand{\ind}[1]{\mathbf1_{\{#1\}}}
\newcommand{\ph}{\phi}
\newcommand{\Span}{\operatorname{span}}

\newtheorem{theorem}{Theorem}[section]
\newtheorem{lemma}[theorem]{Lemma}
\newtheorem{proposition}[theorem]{Proposition}

\theoremstyle{remark}\newtheorem{remark}[theorem]{Remark}

\numberwithin{equation}{section}
\allowdisplaybreaks
\definecolor{accent}{RGB}{23,85,120}
\definecolor{warm}{RGB}{165,70,20}
\hypersetup{pdfauthor={Yi Zheng, Jiansong Deng, Xinyu Wu},pdftitle={Geometry-Dependent Dimension Jumps in Planar Cr Splines at Degree 2r+1}}

\begin{document}

\title{Geometry-Dependent Dimension Jumps in Planar $C^r$ Splines at Degree $2r+1$}
\author{Yi Zheng\\{\small\texttt{zhengyi2023@mail.ustc.edu.cn}}
\and Jiansong Deng\\{\small\texttt{dengjs@ustc.edu.cn}}
\and Xinyu Wu\\{\small\texttt{xinyu97@ustc.edu.cn}}}\date{}
\maketitle
\begin{abstract}
We study geometric dependence of the dimensions of planar spline spaces $S_{2r+1}^r$ through families with a single moving interior vertex and determine their dimensions throughout the admissible parameter domains. For $r=1$, we extend the Wu--Deng one-parameter classification to the interior of the central quadrilateral. The dimension increases by one at exactly two isolated points, one due to a local direction-count change and the other to a global compatibility defect. For $r=2$, subdividing one triangle of the Yuan--Stillman configuration gives a twelve-triangle family with fixed local direction counts; its dimension increases by one precisely on an open line segment. For every $r\ge3$, we construct a family with $2r+12$ triangles and fixed local direction counts, proving geometric dependence at degree $2r+1$. The moving vertex lies in an open triangle, and the dimension exceeds Schumaker's lower bound by one exactly on its symmetry axis. The lower collinear vertices may have arbitrary strictly ordered heights. The low-order proofs use exact cofactor elimination. For higher orders, weighted interpolation and block elimination bound the global defect by one, and explicit extension and obstruction arguments determine the exceptional set.
\end{abstract}
\noindent\textbf{Keywords:} bivariate splines; dimension jump; smoothing cofactors; interpolation; global compatibility.
\par\smallskip\noindent\textbf{Mathematics Subject Classification (2020):} 41A15, 65D07.

\section{Introduction}\label{sec:intro}
Let $\Delta$ be a finite conforming triangulation of a polygonal disk $\Omega\subset\R^2$, with nondegenerate triangles. For integers $d,r\ge0$, write

\[
S_d^r(\Delta)=\{s\in C^r(\Omega):s|_\tau\in\Pa_d\text{ for every }\tau\in\Delta\},
\]
where $\Pa_d$ denotes the real bivariate polynomials of total degree at most $d$. In this paper, geometric dependence with fixed local direction counts means that two nondegenerate geometric realizations of the same abstract triangulation, with the same number of distinct edge directions at each interior vertex, give spline spaces of different dimensions. We study this dependence at degree $d=2r+1$.

Schumaker's lower bound includes a local correction determined by the number of distinct edge directions at each interior vertex~\cite[Theorem~3.1]{Sch84}. For triangulations of a simply connected polygonal domain, equality holds for $d\ge4r+1$ by Alfeld and Schumaker~\cite[Theorem~2.1]{AS87}, and for $d\ge3r+2$ by Hong~\cite[Theorem~3.1 and Corollary~3.1]{Hong}. At $d=3r+1$, Alfeld and Schumaker~\cite[Definition~2.1 and Theorem~2.2]{AS90} prove equality for a dense class satisfying an additional condition on adjacent edge directions. Wu and Deng~\cite[Theorem~1.1 and Corollary~1.2]{Universal} prove attainment for every nondegenerate triangulation of a polygonal disk at $d=3r+1$ for $r\ge3$; together with the low-order cases and Hong\textquoteright s theorem, this gives attainment for all $r\ge1$ and $d\ge3r+1$.

At degree $2r$, geometric dependence is exemplified by the Morgan--Scott triangulation~\cite{MS75}. Diener~\cite[Theorems~1--2]{Diener} establishes dimension instability on this configuration for every $r\ge1$, under the assumption that no two edges incident to an interior vertex are collinear. For these standard configurations, Wu, Zheng, and Deng~\cite[Theorem~2.1 and Corollary~2.2]{WZD} give a geometric criterion: concurrence in the projective plane of three associated lines is sufficient for a dimension increase of one for every $r\ge1$, and is also necessary when $r>2$.

At degree $2r+1$, counterexamples to attainment of Schumaker's lower bound arise from the Schenck--Stiller vanishing conjecture~\cite[Conjecture~2.1]{Stiller}. Yuan and Stillman~\cite[Sections~1.3 and~2.3]{YS} disprove the conjecture with a $C^2$ quintic example. Schenck, Stillman, and Yuan~\cite[Theorem~1.1]{SSY} establish nonvanishing over a larger degree range, excluding a universal threshold $2r+c$ with $c$ independent of $r$. These results establish strict excess over the lower bound on specified triangulations. Geometric dependence asks how this excess varies between realizations with the same local direction counts.

For $r=1,2$, explicit examples already describe such variation. For $C^1$ cubics, Billera~\cite{Billera} proves the generic dimension formula, while Wu and Deng~\cite[Theorem~4.1]{WD} determine the dimension along a horizontal motion of one vertex in an eighteen-triangle complex. Their classification distinguishes a local singularity from a global compatibility defect. For $C^2$ quintics, Wu and Deng~\cite[Example~5.4]{Landscape} exhibit two realizations of a thirteen-triangle complex with identical local direction counts and dimensions $39$ and $40$, obtained by changing three interior vertices. These examples lead to the study of single-vertex families for which the dimension can be determined at every admissible position. A further objective is a construction exhibiting geometric dependence at degree $2r+1$ for every higher smoothness order.

The main result of this paper establishes geometric dependence of the dimension of $S_{2r+1}^r$ for every $r\ge3$. We construct a family of nondegenerate triangulations with fixed combinatorics and local direction counts, in which one interior vertex moves inside a triangle. The dimension is determined throughout this parameter domain: it exceeds Schumaker's lower bound by one on an explicitly specified open line segment and equals the bound elsewhere.

We also obtain full parameter-domain classifications in the two low-order cases. For $r=1$, we extend the horizontal classification of Wu--Deng to the interior of the central quadrilateral and show that the exceptional set consists of exactly two isolated points, with their local and global origins distinguished. For $r=2$, we subdivide one triangle of the Yuan--Stillman configuration by inserting a movable vertex. The resulting twelve-triangle family has fixed local direction counts, and its dimension jumps precisely on an open line segment inside the subdivided triangle.

Our proofs use smoothing cofactor equations~\cite[Sections~1--2]{LW} to separate local vertex contributions from the global defect. This decomposition also appears in the homological approach of Billera~\cite{Billera} and Schenck and Stillman~\cite[Section~3]{SS97}; see~\cite[Sections~2 and~5]{Survey} for an account. Exact matrix elimination gives the low-order classifications. For $r\ge3$, we reduce a global relation to two endpoint components and show that the admissible endpoint space has dimension at most one. An explicit extension on the axis and a nonzero compatibility obstruction off the axis then determine the exceptional set.

Section~\ref{sec:prelim} introduces the cofactor formulation and the global defect. Sections~\ref{sec:cubic}--\ref{sec:high} describe the three families and state their dimension formulas. Section~\ref{sec:discussion} discusses the design of the constructions. Section~\ref{sec:conclusion} summarizes the results and presents further questions. Appendices~\ref{app:low} and~\ref{app:high} contain the proofs for the low-order and higher-order families, respectively.

\section{Smoothing Cofactors and Dimension Decomposition}\label{sec:prelim}
This section expresses the spline dimension in terms of the smoothing cofactor matrix and separates the local vertex contributions from the global defect.

In the dimension formulas and proofs below, $r\ge1$, $d=2r+1$, and $m=r+1$. Let $E^\circ,V^\circ$ be the interior edges and vertices, with cardinalities $f_1^\circ,f_0^\circ$, and let $f_2$ denote the number of triangles. Homogenization identifies $\Pa_j$ with $S_j=\R[x,y,z]_j$, the homogeneous polynomials of degree $j$. Choose a nonzero homogeneous equation $\ell_e$ for each edge and incidence signs $\varepsilon_{ve}$, opposite at the two endpoints of an interior edge. Define
\begin{equation}\label{eq:C}
C_\Delta:\bigoplus_{e\in E^\circ}S_r\longrightarrow\bigoplus_{v\in V^\circ}S_d,
\qquad (C_\Delta q)_v=\sum_{e\ni v}\varepsilon_{ve}\ell_e^m q_e.
\end{equation}For a cofactor vector $q$, we call $(C_\Delta q)_v$ its image component at the vertex $v$.
\begin{proposition}\label{prop:dimension}
For a conforming triangulation of a polygonal disk,
\begin{equation}\label{eq:dimension}
\dim S_{2r+1}^r(\Delta)=\binom{2r+3}{2}+f_1^\circ\binom{r+2}{2}-\rank C_\Delta.
\end{equation}
\end{proposition}
See~\cite[Proposition~2.2 and Appendix~A.2.2]{Universal}.

Let $J_v$ be the image of the vertex block of $C_\Delta$. In coordinates centered at $v$, it is the degree-at-most-$d$ part of the ideal generated by the incident $m$th powers.
After dehomogenizing in coordinates centered at $v$, write $(J_v)_n$ for the binary homogeneous part of total degree $n$. Thus $J_v=\bigoplus_{n=m}^d(J_v)_n$ as a vector space. This local grading is distinct from the common homogenized degree $d$ in~\eqref{eq:C}.
Here $n_v$ counts the distinct unoriented supporting lines through $v$; opposite collinear edges count as one direction.
\begin{lemma}\label{lem:local}
If $n_v$ distinct supporting directions meet at $v$, then
\begin{equation}\label{eq:local}
\dim(J_v)_{m+j}=\min\{n_v(j+1),m+j+1\}\quad(0\le j\le r).
\end{equation}
\end{lemma}
See~\cite[Proposition~2.1 and Appendix~A.1]{Universal}, with the local grading used above.
Set
\begin{gather}
N_r=\binom{2r+3}{2}-\binom{r+2}{2},\qquad
\sigma_v=\sum_{j=1}^{r+1}(r+1-j(n_v-1))_+,\nonumber\\
L_r(\Delta)=\binom{2r+3}{2}+f_1^\circ\binom{r+2}{2}-f_0^\circ N_r+\sum_v\sigma_v,\label{eq:L}\\
\delta_r(\Delta)=\sum_v\dim J_v-\rank C_\Delta\ge0.\label{eq:delta}
\end{gather}
Here $(a)_+=\max(a,0)$. Lemma~\ref{lem:local} gives $\dim J_v=N_r-\sigma_v$, so Proposition~\ref{prop:dimension} yields the exact identity
\begin{equation}\label{eq:split}
\dim S_{2r+1}^r(\Delta)=L_r(\Delta)+\delta_r(\Delta).
\end{equation}
Thus $L_r$ is the Schumaker local lower bound at degree $2r+1$, and $\delta_r$ measures the global defect.

\section{The cubic family}\label{sec:cubic}

Use the fixed vertices of Wu--Deng~\cite[Section~3]{WD}, but allow the central vertex to move in two dimensions. The interior vertices are
\begin{align*}
v_0&=(-3,0),&v_1&=(-3/4,1),&v_2&=(24/55,28/55),\\
v_3&=(3,0),&v_4&=(24/55,-28/55),&v_5&=(-3/4,-1),
\end{align*}
and $v_6=Q=(h,k)$. The boundary vertices are
\[
v_7=(0,4),\ v_8=(-6,4),\ v_9=(-6,-4),\
v_{10}=(0,-4),\ v_{11}=(6,-4),\ v_{12}=(6,4).
\]
The faces, denoted by their vertex indices, are
\begin{equation}\label{eq:cubicfaces}
\begin{gathered}
(0,1,5),(1,2,6),(2,4,6),(2,3,4),(4,5,6),(5,1,6),\\
(7,0,1),(7,1,2),(7,2,3),(10,3,4),(10,4,5),(10,5,0),\\
(0,7,8),(0,8,9),(0,9,10),(3,10,11),(3,11,12),(3,12,7).
\end{gathered}
\end{equation}
Denote this triangulation by $\mathcal T(h,k)$. The admissible set is the interior $\mathcal D_1$ of $[v_1,v_2,v_4,v_5]$. Its four supporting lines are $4h+3=0$, $55h-24=0$, and $12h\pm29k-20=0$. Connecting any interior point to the four vertices gives a nondegenerate fan; the other fourteen triangles remain fixed. Hence the interior of this quadrilateral is the full admissible domain.
\begin{theorem}\label{thm:cubic}
For the triangulation defined above, let
\[
\mathcal D_1=\{(h,k):-3/4<h<24/55,\ |k|<(20-12h)/29\}.
\]
For every $Q=(h,k)$ in this domain,
\[
\dim S_3^1(\mathcal T(h,k))=33+\ind{Q=G\text{ or }Q=O},\qquad
G=(1/5,0),\quad O=(3/83,0).
\]
The local lower bound is $33+\ind{Q=O}$ and the global defect is $\ind{Q=G}$.
\end{theorem}
\begin{remark}
The two exceptional points have different causes. At $O$, the four radial edges lie on two lines, so $n_Q=2$ and the local lower bound increases by one. At $G$ and in a sufficiently small neighborhood, the direction counts at $(v_0,\ldots,v_6)$ remain $(4,5,5,4,5,5,4)$; only the global defect changes. Thus $G$ exhibits geometric dependence with fixed local direction counts, whereas $O$ is a local singularity.
\end{remark}

\begin{inplacefigure}\centering
\begin{minipage}[c]{.61\linewidth}\centering
\begin{minipage}[c][5.0cm][c]{\linewidth}\centering
\begin{tikzpicture}[scale=.55,font=\small]
\foreach \i/\x/\y in {0/-3/0,1/-.75/1,2/.436364/.509091,3/3/0,4/.436364/-.509091,5/-.75/-1,6/.2/0,7/0/4,8/-6/4,9/-6/-4,10/0/-4,11/6/-4,12/6/4}
{\coordinate(v\i)at(\x,\y);}
\fill[accent!12](v1)--(v2)--(v4)--(v5)--cycle;
\foreach \a/\b in {0/1,0/5,1/5,1/2,2/4,2/3,3/4,4/5,0/7,1/7,2/7,3/7,3/10,4/10,5/10,0/10,7/8,8/9,9/10,10/11,11/12,12/7,0/8,0/9,3/11,3/12}
{\draw[gray](v\a)--(v\b);}
\foreach \i in {1,2,4,5}{\draw[warm,thick](v\i)--(v6);}
\foreach \i in {0,...,12}{\fill(v\i)circle(1.6pt);}
\node[left]at(v0){$v_0$};\node[right]at(v3){$v_3$};
\node[above left]at(v1){$v_1$};\node[above right]at(v2){$v_2$};
\node[below right]at(v4){$v_4$};\node[below left]at(v5){$v_5$};
\draw[warm](v6)--(1.3,-1.6)node[right]{$Q$};
\node[above]at(v7){$v_7$};\node[below]at(v10){$v_{10}$};
\node[left]at(v8){$v_8$};\node[left]at(v9){$v_9$};
\node[right]at(v11){$v_{11}$};\node[right]at(v12){$v_{12}$};
\end{tikzpicture}
\end{minipage}\par\smallskip
(a)
\end{minipage}\hfill
\begin{minipage}[c]{.37\linewidth}\centering
\begin{minipage}[c][5.0cm][c]{\linewidth}\centering
\begin{tikzpicture}[xscale=2.65,yscale=1.8,font=\small]
\fill[accent!9](-.75,1)--(.436364,.509091)--(.436364,-.509091)--(-.75,-1)--cycle;
\draw(-.75,1)--(.436364,.509091)--(.436364,-.509091)--(-.75,-1)--cycle;
\draw[gray,dashed](-.75,1)--(.436364,-.509091);
\draw[gray,dashed](-.75,-1)--(.436364,.509091);
\draw[->](-.85,0)--(.62,0)node[right]{$h$};
\node[circle,fill=warm,inner sep=0pt,minimum size=3.6pt]at(.2,0){};\node[circle,fill=accent,inner sep=0pt,minimum size=3.6pt]at(.036145,0){};
\draw[warm](.2,0)--(.2,.85)node[above]{$G$};
\draw[accent](.036145,0)--(-.15,-.8)node[below]{$O$};
\node at(-.47,.48){$\mathcal D_1$};
\end{tikzpicture}
\end{minipage}\par\smallskip
(b)
\end{minipage}
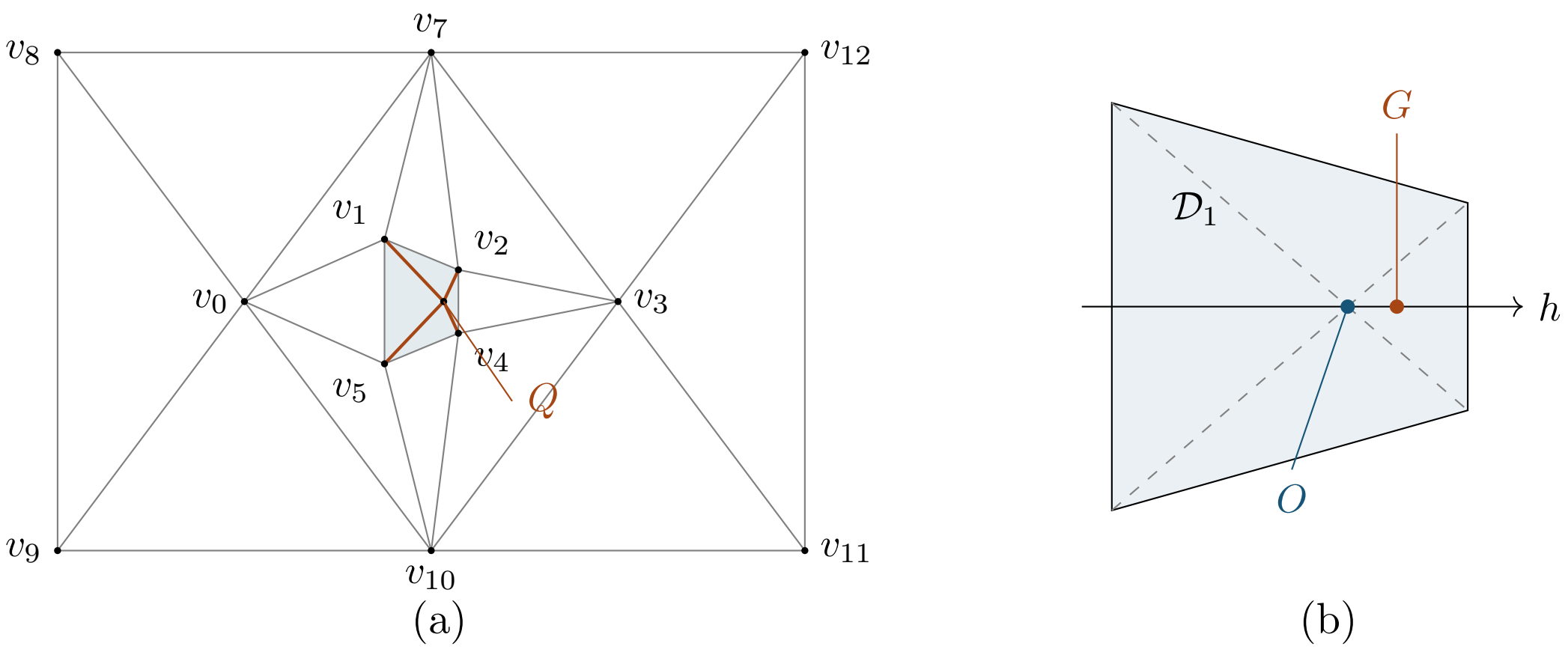
\caption{(a) The cubic mesh at $Q=G$. (b) Its admissible quadrilateral. The global defect occurs only at $G$; the diagonal intersection $O$ contributes a local correction.}\label{fig:cubic}
\end{inplacefigure}

The proof is given in Appendix~\ref{app:cubicproof}.

\section{The quintic family}\label{sec:quintic}

Fix
\begin{align*}
v_0&=(0,-1),&v_1&=(-2,0),&v_2&=(2,0),\\
v_3&=(-4,2),&v_4&=(-4,-2),&v_5&=(0,-2),\\
v_6&=(4,-2),&v_7&=(4,2),&v_8&=(0,2),
\end{align*}
and put $v_9=Q=(h,k)$. The twelve faces are
\begin{equation}\label{eq:quinticfaces}
\begin{gathered}
(3,4,1),(3,8,1),(4,5,1),(1,0,5),\\
(2,7,6),(2,8,7),(2,6,5),(2,0,8),(2,5,0),\\
(1,8,9),(8,0,9),(0,1,9).
\end{gathered}
\end{equation}
The last three subdivide $[v_1,v_8,v_0]$. The barycentric coordinates of $Q$ in this triangle are
\[
\left(-\frac h2,\frac{k+1+h/2}{3},\frac{2+h-k}{3}\right).
\]
These three coordinates are strictly positive precisely when $Q$ belongs to the admissible domain
\[
\mathcal D_2=\{(h,k):-2<h<0,\ -1-h/2<k<h+2\}.
\]
For every $(h,k)\in\mathcal D_2$, the fixed outer triangles and this internal subdivision give a conforming nondegenerate triangulation $\Delta_2(h,k)$.
Throughout $\mathcal D_2$, the direction counts at $(v_0,v_1,v_2,Q)$ are $(4,4,3,3)$. Thus the local lower bound is constant as $Q$ moves.
\begin{theorem}\label{thm:quintic}
For every $(h,k)\in\mathcal D_2$, the triangulation defined above satisfies
\[
\dim S_5^2(\Delta_2(h,k))=53+\ind{13h-5k+10=0}.
\]
\end{theorem}

\begin{inplacefigure}\centering
\begin{minipage}[c]{.60\linewidth}\centering
\begin{minipage}[c][4.9cm][c]{\linewidth}\centering
\begin{tikzpicture}[scale=.67,font=\small]
\foreach \i/\x/\y in {0/0/-1,1/-2/0,2/2/0,3/-4/2,4/-4/-2,5/0/-2,6/4/-2,7/4/2,8/0/2,9/-.5/.7}
{\coordinate(v\i)at(\x,\y);}
\fill[accent!12](v1)--(v8)--(v0)--cycle;
\foreach \a/\b in {3/4,3/8,4/5,5/6,6/7,7/8,3/1,4/1,5/1,8/1,0/1,0/5,0/8,0/2,2/5,2/6,2/7,2/8}
{\draw[gray](v\a)--(v\b);}
\foreach \i in {0,1,8}{\draw[warm,thick](v\i)--(v9);}
\foreach \i in {0,...,9}{\fill(v\i)circle(1.6pt);}
\node[below left]at(v0){$v_0$};\node[left]at(v1){$v_1$};\node[right]at(v2){$v_2$};
\node[above]at(v8){$v_8$};\node[right,warm]at(v9){$Q$};
\node[above left]at(v3){$v_3$};\node[below left]at(v4){$v_4$};\node[below]at(v5){$v_5$};
\node[below right]at(v6){$v_6$};\node[above right]at(v7){$v_7$};
\end{tikzpicture}
\end{minipage}\par\smallskip
(a)
\end{minipage}\hfill
\begin{minipage}[c]{.38\linewidth}\centering
\begin{minipage}[c][4.9cm][c]{\linewidth}\centering
\begin{tikzpicture}[scale=1.0,font=\small]
\fill[accent!9](-2,0)--(0,2)--(0,-1)--cycle;
\draw(-2,0)--(0,2)--(0,-1)--cycle;
\draw[warm,thick](-.967742,-.516129)--(0,2);
\draw[warm,fill=white](-.967742,-.516129)circle(1.8pt);
\draw[warm,fill=white](0,2)circle(1.8pt);
\node[right]at(0,2){$(0,2)$};
\node[below left=2pt,font=\small,inner sep=1pt]at(-.967742,-.516129){$(-\frac{30}{31},-\frac{16}{31})$};
\node[left]at(-1.1,.6){$\mathcal D_2$};
\node[right,warm,font=\small,inner sep=1pt]at(.12,.75){$13h-5k+10=0$};
\end{tikzpicture}
\end{minipage}\par\smallskip
(b)
\end{minipage}
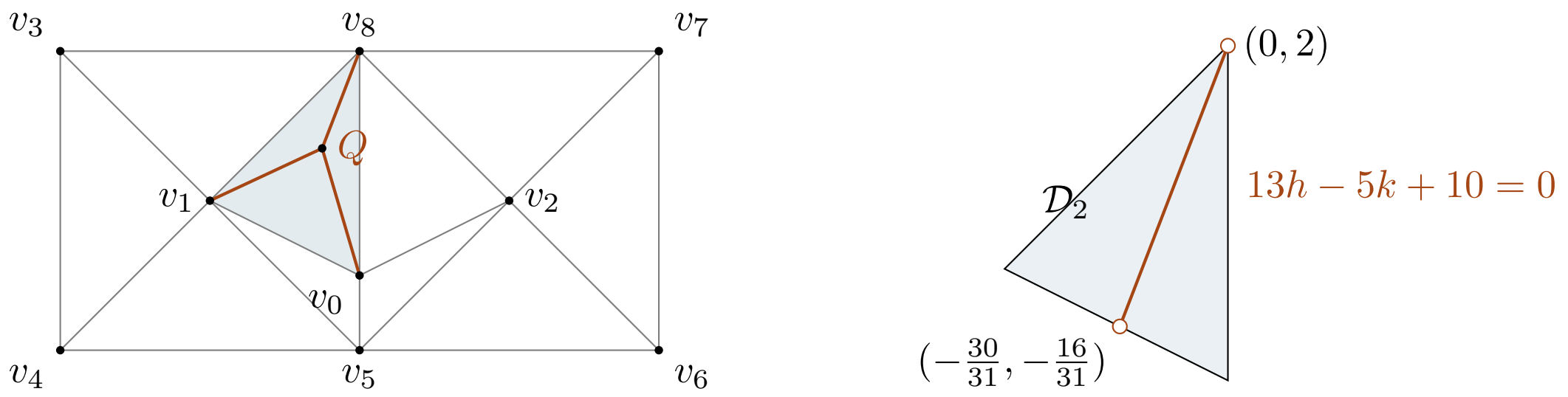
\caption{(a) The twelve-triangle quintic mesh. (b) The admissible domain and exceptional open segment. The endpoints are excluded from the admissible domain.}\label{fig:quintic}
\end{inplacefigure}

The proof is given in Appendix~\ref{app:quinticproof}.

\section{The higher-order family}\label{sec:high}

Fix $r\ge3$. The boundary vertices of $[-3,3]\times[-2,2]$, in counterclockwise order, are
\begin{gather*}
v_0=(-3,-2),\ v_1=(-1,-2),\ v_2=(1,-2),\ v_3=(3,-2),\\
v_4=(3,2),\ v_5=(1,2),\ v_6=(-1,2),\ v_7=(-3,2).
\end{gather*}
Add interior vertices
\[
A=v_8=(-1,0),\quad B=v_9=(1,0),\quad X=v_{10}=(0,-1),\quad Y=v_{11}=(0,1).
\]
Before inserting any of the vertices $P_i$ or $Q$, define the initial triangulation $\Delta_{\rm old}$ by the fourteen faces
\begin{equation}\label{eq:oldfaces}
\begin{gathered}
(8,1,0),(8,0,7),(8,7,6),(9,3,2),(9,4,3),(9,5,4),\\
(10,1,2),(10,2,9),(10,9,8),(10,8,1),\\
(11,8,9),(11,9,5),(11,5,6),(11,6,8).
\end{gathered}
\end{equation}
Its only interior vertices are $A,B,X,Y$; see Figure~\ref{fig:high}(a). An old edge means an edge of $\Delta_{\rm old}$, all of which are retained in the subsequent insertions. An old direction at a vertex is a supporting-line direction of an incident edge in this initial triangulation, with collinear edges counted only once. Thus $A$ and $B$ each have four old directions (horizontal, vertical, and slopes $1$ and $-1$), whereas $X$ and $Y$ each have two (slopes $1$ and $-1$).

Choose lower heights satisfying
\begin{equation}\label{eq:heights}
-1<t_0<\cdots<t_{r-3}<0,\qquad P_i=(0,t_i)\quad(0\le i\le r-3).
\end{equation}
Insert $P_0$ into $[A,B,X]$, joining it to the three vertices, and successively insert $P_i$ into $[A,B,P_{i-1}]$. Finally insert $Q=(h,k)$ into $[A,B,Y]$. Fixing the ordered lower heights, denote the resulting triangulation by $\Delta_r(h,k)$. Put $t_{-1}=-1$, $P_{-1}=X$. The barycentric coordinates for the $i$th lower insertion are
\[
\left(\frac{1-t_i/t_{i-1}}2,\frac{1-t_i/t_{i-1}}2,\frac{t_i}{t_{i-1}}\right).
\]
They are positive because $0<t_i/t_{i-1}<1$. The upper insertion has coordinates $((1-k-h)/2,(1-k+h)/2,k)$, giving the admissible domain below.
Throughout the admissible domain, $n_A=n_B=r+3$, and each of the other $r+1$ interior vertices has three distinct directions. These counts are independent of $Q$.

\begin{inplacefigure}\centering
\begin{minipage}[c]{.57\linewidth}\centering\begin{minipage}[c][4cm][c]{\linewidth}\centering
\begin{tikzpicture}[scale=.72,font=\small]
\foreach \i/\x/\y in {0/-3/-2,1/-1/-2,2/1/-2,3/3/-2,4/3/2,5/1/2,6/-1/2,7/-3/2,8/-1/0,9/1/0,10/0/-1,11/0/1}
{\coordinate(v\i)at(\x,\y);}
\fill[accent!10](v8)--(v9)--(v11)--cycle;
\foreach \a/\b in {0/1,1/2,2/3,3/4,4/5,5/6,6/7,7/0,8/1,8/0,8/7,8/6,9/3,9/2,9/4,9/5,10/1,10/2,10/9,10/8,8/9,11/8,11/9,11/5,11/6}
{\draw[gray](v\a)--(v\b);}
\foreach \i in {0,...,11}{\fill(v\i)circle(1.4pt);}

\node[left]at(v8){$A$};\node[right]at(v9){$B$};
\node[below]at(v10){$X$};\node[above]at(v11){$Y$};

\end{tikzpicture}
\end{minipage}\par\smallskip(a) $\Delta_{\rm old}$
\end{minipage}\hfill
\begin{minipage}[c]{.40\linewidth}\centering\begin{minipage}[c][4cm][c]{\linewidth}\centering
\begin{tikzpicture}[scale=1.35,font=\small]
\coordinate(A)at(-1,0);\coordinate(B)at(1,0);\coordinate(X)at(0,-1);\coordinate(Y)at(0,1);
\coordinate(p0)at(0,-.5);\coordinate(p1)at(0,-.166667);\coordinate(q)at(.25,.55);
\fill[accent!9](A)--(B)--(Y)--cycle;
\draw[gray](A)--(X)--(B)--(Y)--cycle;\draw[gray](A)--(B);
\foreach \p in {p0,p1}{\draw[accent](A)--(\p)--(B);\node[circle,fill=accent,inner sep=0pt,minimum size=3pt]at(\p){};}
\draw[accent](X)--(p1);
\draw[warm](A)--(q)--(B);\draw[warm](q)--(Y);
\foreach \p in {A,B,X,Y}{\node[circle,fill,inner sep=0pt,minimum size=3pt]at(\p){};}
\node[left]at(A){$A$};\node[right]at(B){$B$};\node[below]at(X){$X$};\node[above]at(Y){$Y$};
\node[circle,fill=warm,inner sep=0pt,minimum size=3pt]at(q){};\node[right,warm]at(q){$Q$};
\node[below right,accent,inner sep=1.5pt]at(p0){$P_0$};\node[below right,accent,inner sep=1.5pt]at(p1){$P_1$};

\end{tikzpicture}
\end{minipage}\par\smallskip(b) $r=4$
\end{minipage}
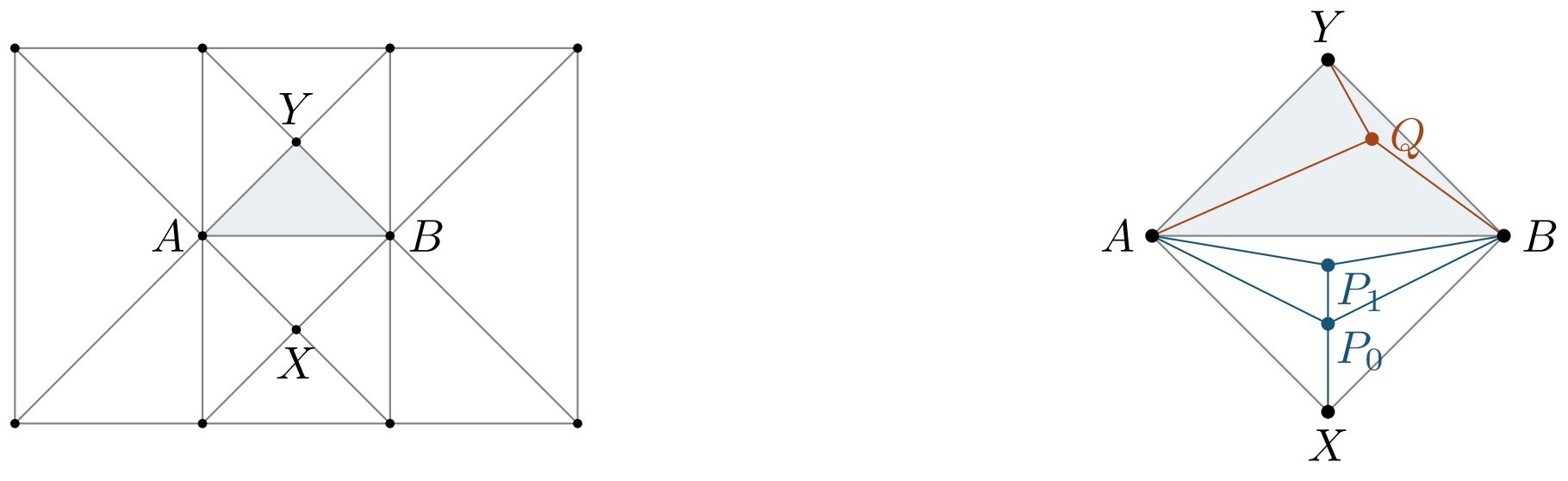
\caption{(a) The initial fourteen-triangle mesh $\Delta_{\rm old}$. (b) The central region $AXBY$ after refinement for $r=4$. Gray edges belong to $\Delta_{\rm old}$, blue edges arise from the lower insertions, and orange edges join $Q$ to $A,B,Y$. Here $t_0=-1/2$, $t_1=-1/6$; the heights may be chosen as in~\eqref{eq:heights}.}\label{fig:high}
\end{inplacefigure}

\begin{theorem}\label{thm:high}
Let $r\ge3$ and fix any ordered heights satisfying~\eqref{eq:heights}.
For the triangulation defined above, let
\[
\mathcal D=\{(h,k):0<k<1,\ |h|<1-k\}.
\]
For every $Q=(h,k)$ in this domain,
\[
\dim S_{2r+1}^r(\Delta_r(h,k))=
11\binom{r+2}{2}+(r+1)\left\lfloor\frac{r^2}{4}\right\rfloor+\ind{h=0}.
\]
\end{theorem}
\subsection*{Proof sketch}
Appendix~\ref{app:bound} first removes the local relations and reduces a global left relation to its components at $A,B$. The lower-chain conditions are then written as univariate interpolation conditions; weighted interpolation and block elimination give an endpoint space of dimension at most one. The injectivity of the endpoint restriction yields $\delta_r\le1$. Appendix~\ref{sec:obstruction} constructs, for $h=0$, an explicit left-kernel vector whose class modulo local relations is nonzero, proving $\delta_r\ge1$. Appendix~\ref{app:offaxis} pairs the candidate endpoint relation, unique up to a scalar, with the image of a constructed cofactor vector. The resulting compatibility obstruction is nonzero when $h\ne0$, so no nonzero global class remains and $\delta_r=0$.

\begin{inplacefigure}\centering
\begin{minipage}[t]{.32\linewidth}\centering
\begin{tikzpicture}[scale=1.65,font=\small]
\path[use as bounding box](-1.25,-.45)rectangle(1.25,1.3);
\fill[accent!18](-1,0)--(1,0)--(0,1)--cycle;
\node[accent]at(0,.35){$\mathcal D$};
\draw[gray](-1,0)--(1,0)--(0,1)--cycle;
\node[left]at(-1,0){$A$};\node[right]at(1,0){$B$};\node[above]at(0,1){$Y$};
\end{tikzpicture}\par
(a) Appendix~\ref{app:bound}\par $\delta_r\le1$
\end{minipage}\hfill
\begin{minipage}[t]{.32\linewidth}\centering
\begin{tikzpicture}[scale=1.65,font=\small]
\path[use as bounding box](-1.25,-.45)rectangle(1.25,1.3);
\draw[gray](-1,0)--(1,0)--(0,1)--cycle;
\node[left]at(-1,0){$A$};\node[right]at(1,0){$B$};\node[above]at(0,1){$Y$};
\draw[warm,very thick](0,0)--(0,1);
\draw[warm,fill=white](0,0)circle(1.5pt);\draw[warm,fill=white](0,1)circle(1.5pt);
\node[right,warm]at(0,.45){$h=0$};
\end{tikzpicture}\par
(b) Appendix~\ref{sec:obstruction}\par $\delta_r\ge1$
\end{minipage}\hfill
\begin{minipage}[t]{.32\linewidth}\centering
\begin{tikzpicture}[scale=1.65,font=\small]
\path[use as bounding box](-1.25,-.45)rectangle(1.25,1.3);
\fill[accent!18](-1,0)--(0,0)--(0,1)--cycle;
\fill[accent!18](1,0)--(0,0)--(0,1)--cycle;
\draw[white,line width=2.5pt](0,0)--(0,1);
\draw[gray,dashed](0,0)--(0,1);
\draw[gray](-1,0)--(1,0)--(0,1)--cycle;
\node[left]at(-1,0){$A$};\node[right]at(1,0){$B$};\node[above]at(0,1){$Y$};
\node at(0,-.3){$h\ne0$};
\end{tikzpicture}\par
(c) Appendix~\ref{app:offaxis}\par $\delta_r=0$
\end{minipage}
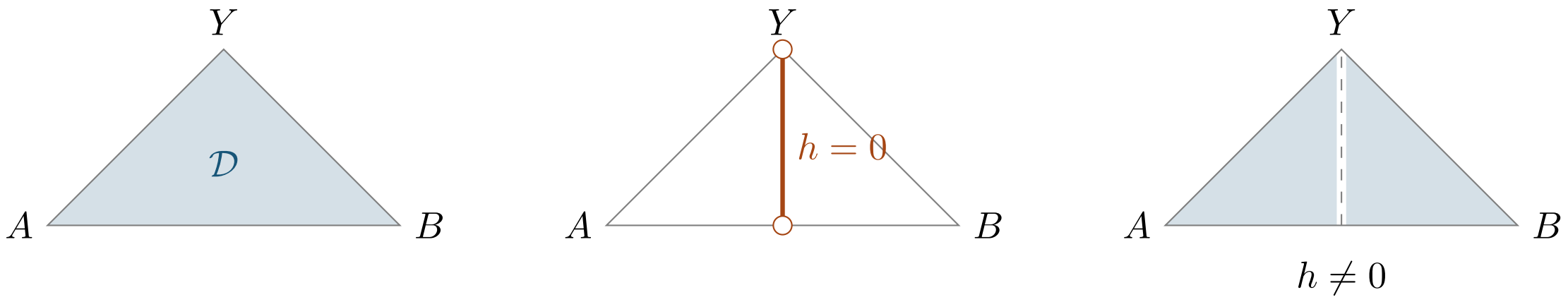
\caption{The parameter regions in the three steps of the proof. (a) The bound $\delta_r\le1$ holds throughout the open triangle $\mathcal D$. (b) A nonlocal left relation exists on the open segment $h=0$. (c) The defect vanishes in the two regions $h\ne0$.}\label{fig:proofroute}
\end{inplacefigure}

The proof is given in Appendix~\ref{app:high}.

\section{Discussion}\label{sec:discussion}
The constructions start from a common question: how can one vary a global compatibility condition by moving a single vertex while keeping the local direction counts fixed? The low-order examples suggest local refinement, and the higher-order family makes this refinement depend on the smoothness order.
\subsection{Local refinement of known examples}
For $r=1$, we retain the Wu--Deng configuration~\cite{WD} and extend its horizontal one-parameter classification to the whole central quadrilateral. This determines whether the known exceptional positions are isolated or belong to larger exceptional sets. Keeping the surrounding mesh fixed and moving a single interior vertex makes it possible to describe the exceptional set in a two-dimensional parameter domain.

For $r=2$, we start with the ten-triangle example of Yuan--Stillman~\cite[Figure~1]{YS}, redrawn in Figure~\ref{fig:yuanmotion}. We first try moving $v_0$ vertically while keeping the other vertices fixed. With the interior direction counts unchanged, this motion does not produce a dimension jump. We therefore subdivide a triangle and introduce a new moving vertex. The surrounding mesh stays fixed, while the three new edges provide compatibility conditions that vary with the inserted point.

\begin{inplacefigure}\centering
\begin{tikzpicture}[scale=.65,font=\small]
\foreach \i/\x/\y in {0/0/-1,1/-2/0,2/2/0,3/-4/2,4/-4/-2,5/0/-2,6/4/-2,7/4/2,8/0/2}{\coordinate(v\i)at(\x,\y);}
\fill[accent!9](v1)--(v8)--(v0)--cycle;
\foreach \a/\b in {3/4,3/8,4/5,5/6,6/7,7/8,3/1,4/1,5/1,8/1,0/1,0/5,0/8,0/2,2/5,2/6,2/7,2/8}{\draw[gray](v\a)--(v\b);}
\foreach \i in {0,...,8}{\fill(v\i)circle(1.5pt);}
\node[below left]at(v0){$v_0=(0,-1)$};\node[left]at(v1){$v_1$};\node[right]at(v2){$v_2$};
\node[above left]at(v3){$v_3$};\node[below left]at(v4){$v_4$};\node[below]at(v5){$v_5$};\node[below right]at(v6){$v_6$};\node[above right]at(v7){$v_7$};\node[above]at(v8){$v_8$};
\draw[warm,<->,thick](.16,-1.65)--(.16,.25);\node[right,warm]at(.25,-.55){$(0,\kappa)$};
\end{tikzpicture}
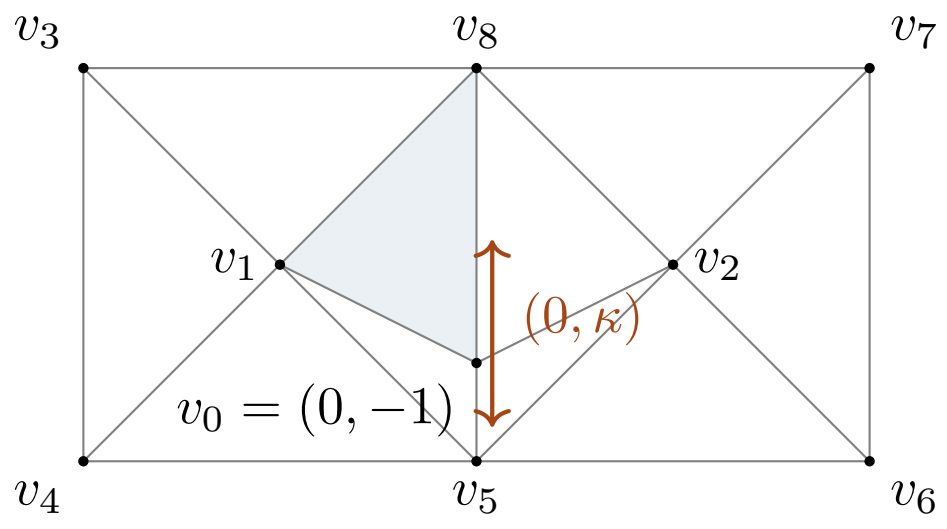
\caption{The Yuan--Stillman triangulation, redrawn from~\cite[Figure~1]{YS}. The arrow indicates the vertical motion considered here; the shaded face $[v_1,v_8,v_0]$ is the one subdivided in Section~\ref{sec:quintic}.}\label{fig:yuanmotion}
\end{inplacefigure}

The choice of the triangle is part of the construction. In the face $[v_1,v_8,v_0]$ used here, the interior direction counts remain $(4,4,3,3)$ throughout the motion of $Q$. The dimension formula then identifies the exceptional line $13h-5k+10=0$. This suggests looking for faces whose subdivision preserves the local direction counts and changes the compatibility with the surrounding mesh. Both properties matter: the location of a face in the interior is only an initial geometric criterion, and the effect of its subdivision must be determined from the resulting compatibility conditions.

\subsection{Fixed and moving parts of the higher-order construction}
For $r\ge3$, the construction retains a single moving vertex and allows the fixed part to grow with the smoothness order. Two triangles sharing $AB$ provide a convenient layout: the upper triangle contains the moving point, while the lower triangle admits repeated subdivision. Connecting each inserted point to $A$ and $B$ preserves the same local pattern as the mesh grows. The edge $AB$ separates the two regions geometrically, while smoothness at the common vertices couples their compatibility conditions. This makes it possible to adjust the fixed lower part and the moving upper part separately within one family.

In this layout, the lower part can be used to restrict the possible global relations, and the position of the upper vertex can be used to test their compatibility with the remaining triangles. The proof confirms this role for the present construction: $r-2$ lower insertions restrict the normalized endpoint pairs to a space of dimension at most one; whether a pair extends to the whole triangulation depends on $Q$. The lower heights may be chosen arbitrarily in strict order. Thus the construction provides a way to control the number of possible relations with a fixed part of the mesh while keeping the geometric parameter confined to one triangle.

\clearpage
\section{Conclusion}\label{sec:conclusion}
For every integer $r\ge1$, there exist two nondegenerate realizations of the same abstract triangulation, with the same number of distinct edge directions at each interior vertex, for which $S_{2r+1}^r$ has different dimensions. Thus geometric dependence occurs at degree $2r+1$ for every positive smoothness order. For $r=1$, this follows by comparing $G$ with nearby positions; the exceptional point $O$ is accounted for separately by the change in the local lower bound. For $r=2$, the exceptional set is an open line segment. For every $r\ge3$, the constructed family has dimension one above the local lower bound precisely when the moving vertex lies on the symmetry axis, and attains that bound elsewhere. In each case, the dimension formula covers the full admissible domain of the moving vertex.

\subsection*{Future work}
Increasing smoothness on a fixed mesh can remove its dimension jumps while preserving a positive global defect. For the $r=3$ mesh of Section~\ref{sec:high}, with $P_0=(0,-1/2)$ fixed, exact symbolic elimination gives $\dim S_9^4(\Delta_3(h,k))=178$ throughout $0<k<1$, $|h|<1-k$. The Schumaker lower bound is $176$, so the global defect is identically two. This example raises a general question: for a mesh constructed at smoothness $s$, does the dimension of $S_{2r+1}^r$ remain constant throughout the admissible domain whenever $r>s$? Further questions concern the smallest number of insertions needed for geometric dependence and constructions with larger jumps or more complicated exceptional sets.

\section*{Acknowledgements}
\begingroup\small\emergencystretch=1em
This work was supported by the National Key R\&D Program of China (Nos. 2022YFA1005200 and 2022YFA1005201), the NSF of China (No. 12171453), and the Major Project of Science and Technology Innovation Tackling Plan of Anhui Province (No. 202423e09050003).
\par\endgroup

The authors thank Bingru Huang and Beihui Yuan for helpful discussions.

\section*{Declarations}
The authors declare that they have no conflict of interest.

\appendix

\clearpage
\section{Proofs for the low-order families}\label{app:low}
We specify a common ordering of the matrix entries and a block reduction that will be used for both low-order families.

For the low-order calculations, order interior vertices and edges lexicographically. Row and column sets used to specify matrix minors are indexed from zero. For $e=(i,j)$ with $i<j$, put $A_e=y_i-y_j$, $B_e=x_j-x_i$, and take sign $+1$ at $i$ and $-1$ at $j$. In coordinates $X=x-x_v$, $Y=y-y_v$, the rows at $v$ correspond to the coefficients of $X^\alpha Y^\beta$, ordered by increasing total degree $m,\ldots,d$, and within each degree by decreasing $\alpha$. Cofactor columns are $x^ay^b$ in the analogous order of degrees $0,\ldots,r$. Here $[X^\alpha Y^\beta]$ denotes coefficient extraction. The entries are
\begin{equation}\label{eq:entry}
\varepsilon_{ve}[X^\alpha Y^\beta](A_eX+B_eY)^m(x_v+X)^a(y_v+Y)^b.
\end{equation}
Deleting locally zero rows and changing to centered coordinates preserves rank. Suppose old edges and new radial edges give
\[
C=\begin{pmatrix}C_0&B\\0&D\end{pmatrix},
\]
and let the rows of $U$ span $\ker C_0^{\mathsf T}$. If the columns of $N$ form a basis of $\ker D$, then
\begin{equation}\label{eq:reduction}
\begin{split}
\rank C&=\rank C_0+\rank\begin{pmatrix}UB\\D\end{pmatrix}\\
&=\rank C_0+\rank D+\rank(UBN).
\end{split}
\end{equation}
Indeed, $\ker U=\im C_0$, so the first equality takes the quotient by the old image. For the second, split the new input into $\ker D$ and a complement. The complement contributes $\rank D$ independent new-vertex components; the kernel contributes precisely $\im UBN$ in the old-image quotient.

\subsection{Proof for $r=1$}\label{app:cubicproof}

There are $18$ triangles, $24$ interior edges, and $7$ interior vertices. Each fixed interior vertex has valence at least five, hence at least three directions and local rank seven. At $Q$, two directions occur precisely when the opposite radial edges are collinear in pairs. The relevant diagonals are
\[
332h+261k-12=0,\qquad332h-261k-12=0.
\]
Their intersection is $O=(3/83,0)$. At $O$, the degree-two and degree-three local ranks are two and four; elsewhere they are three and four. Consequently
\begin{equation}\label{eq:cubiclocal}
\dim S_3^1=82-\rank C,\qquad L_1=33+\ind{Q=O}.
\end{equation}

Delete the four radial edges and the rows at $Q$, and denote the resulting $42\times60$ matrix by $C_0$.
\subsubsection{The constant matrix}\label{app:cubicmatrix}

The $42$ entries of the left-kernel row $U$ are grouped into six seven-entry blocks:
\begin{align*}
U_0={}&(-1/6,0,8/27,0,0,0,0),\\
U_1={}&(19/48,11/12,1,-9/64,-19/48,-11/12,-1),\\
U_2={}&(-18272/9075,75776/9075,-14336/3025,\\
&\hspace{5mm}-31584/166375,511616/499125,-2121728/499125,401408/166375),\\
U_3={}&(70/81,0,-1120/729,0,0,0,0),\\
U_4={}&(-18272/9075,-75776/9075,-14336/3025,\\
&\hspace{5mm}-31584/166375,-511616/499125,-2121728/499125,-401408/166375),\\
U_5={}&(19/48,-11/12,1,-9/64,19/48,-11/12,1).
\end{align*}
Expansion using~\eqref{eq:entry} verifies $UC_0=0$. Since $U\ne0$, $\rank C_0\le41$. Take rows $0,\ldots,40$ and columns
\[
\{0,\ldots,11\}\cup\{18,\ldots,44\}\cup\{51,52\}.
\]
The absolute value of this $41\times41$ minor is
\[
\frac{2^{93}3^{70}7^8\,29^3\,47^5}{5^{24}11^{28}}\ne0.
\]
This proves $\rank C_0=41$ and that $U$ spans its full left kernel.

\subsubsection{Two successive compatibility conditions}
After eliminating the fixed columns, we determine the remaining compatibility conditions in degrees two and three.

The matrix $C_0$ and its left-kernel row $U$ are as above. Formula~\eqref{eq:reduction} gives
\[
\rank C=41+\rank E,\qquad E=\begin{pmatrix}g\\D\end{pmatrix},\quad g=UB.
\]
Replace each cofactor $a+bx+cy$ by $a'+bX+cY$, where $a'=a+hb+kc$, $X=x-h$, $Y=y-k$. The change is invertible throughout the domain. Grouping the four constant columns and eight linear columns gives
\begin{equation}\label{eq:cubicblocks}
E\sim\begin{pmatrix}g_2&g_3\\D_2&0\\0&D_3\end{pmatrix},\qquad
D_2\in\R^{3\times4},\quad D_3\in\R^{4\times8}.
\end{equation}
For $i=1,2,4,5$, set $A_i=y_i-k$, $B_i=h-x_i$. The corresponding column of $D_2$ and pair of columns of $D_3$ are
\[
-\begin{pmatrix}A_i^2\\2A_iB_i\\B_i^2\end{pmatrix},\qquad
-\begin{pmatrix}A_i^2&0\\2A_iB_i&A_i^2\\B_i^2&2A_iB_i\\0&B_i^2\end{pmatrix}.
\]
\begin{lemma}\label{lem:cubiccalc}
With the preceding normalization,
\begin{equation}\label{eq:cubich}
\det\begin{pmatrix}g_2\\D_2\end{pmatrix}
=-\frac{10971}{9150625}k(4h+3)(55h-24)(12h-29k-20)(12h+29k-20).
\end{equation}
If $P_3$ consists of the first four columns of $D_3$, then
\begin{equation}\label{eq:cubicpivot}
\det P_3=\frac{6561}{2342560000}(12h+29k-20)^4\ne0\quad\text{on }\mathcal D_1.
\end{equation}
Writing $g_3=(g_{3,L},g_{3,R})$ with four entries in each block, the last four entries of $g_3-g_{3,L}P_3^{-1}D_3$ at $k=0$ are
\begin{align}
z_1&=-\frac{6272(4h+3)(5h-1)(55h-24)^2}{735075(3h-5)},\nonumber\\
z_2&=\frac{12544(5h-1)(55h-24)(83h-3)}{735075(3h-5)},\nonumber\\
z_3&=-\frac{56(4h+3)^2(5h-1)(55h-24)}{243(3h-5)},\nonumber\\
z_4&=\frac{112(4h+3)(5h-1)(83h-3)}{243(3h-5)}.\label{eq:cubicres}
\end{align}
\end{lemma}
\label{app:cubiccalc}
\begin{proof}[Proof of Lemma~\ref{lem:cubiccalc}]
For $i\in\{1,2,4,5\}$ put $A_i=y_i-k$, $B_i=h-x_i$, and $w_i=(A_i^2,2A_iB_i,B_i^2)$. Split the seven entries of the corresponding vertex block of $U$ into $u_i^{(2)}\in\R^3$ and $u_i^{(3)}\in\R^4$, ordered as in~\eqref{eq:entry}. The three centered cofactor columns $1,x-h,y-k$ give the entries
\begin{align*}
g_{2,i}&=u_i^{(2)}\cdot w_i,\\
g_{3,i,x}&=(x_i-h)g_{2,i}
+u_i^{(3)}\cdot(A_i^2,2A_iB_i,B_i^2,0),\\
g_{3,i,y}&=(y_i-k)g_{2,i}
+u_i^{(3)}\cdot(0,A_i^2,2A_iB_i,B_i^2).
\end{align*}
Indeed, at $v_i$ the radial line is $A_iX+B_iY$, and $x-h=(x_i-h)+X$, $y-k=(y_i-k)+Y$. Multiplication by its square separates exactly into the displayed quadratic and cubic terms. All four old-endpoint incidence signs are positive.

Expand the determinant along its first row:
\[
\det\begin{pmatrix}g_2\\D_2\end{pmatrix}
=\sum_{j=1}^4(-1)^{1+j}g_{2,i_j}\det D_2^{\widehat j},
\]
where $(i_1,i_2,i_3,i_4)=(1,2,4,5)$ lists the radial vertices and $D_2^{\widehat j}$ deletes the $j$th column. Substitution of the three-term dot products above and collection of factors gives~\eqref{eq:cubich}. For the pivot matrix, put $\ell_i=A_iX+B_iY$ for $i=1,2$. Its columns are the coefficients of $-X\ell_1^2,-Y\ell_1^2,-X\ell_2^2,-Y\ell_2^2$. Their determinant is
\[
\det P_3=(A_1B_2-A_2B_1)^4.
\]
This identity follows by expanding the displayed four-by-four matrix in~\eqref{eq:cubicblocks}; substitution of the vertex coordinates gives~\eqref{eq:cubicpivot}. To compute the remaining entries, write $D_3=(P_3\mid R_3)$ and denote column $j$ of $R_3$ by $R_{3,j}$. Then
\[
z_j=\left.
\frac{g_{3,R,j}\det P_3-g_{3,L}\operatorname{adj}(P_3)R_{3,j}}
{\det P_3}\right|_{k=0},\qquad 1\le j\le4.
\]
Using the preceding dot products in this formula gives the four factorizations~\eqref{eq:cubicres}. The adjugate identity $P_3\operatorname{adj}(P_3)=(\det P_3)I$ proves that these are exactly the residuals of the row elimination. Division is valid throughout $\mathcal D_1$ by~\eqref{eq:cubicpivot}.

\end{proof}
\begin{proof}[Proof of Theorem~\ref{thm:cubic}]
The two column groups in~\eqref{eq:cubicblocks} are disjoint, so $g\in\row D$ precisely when $g_2\in\row D_2$ and $g_3\in\row D_3$. Every factor in~\eqref{eq:cubich} other than $k$ is nonzero in $\mathcal D_1$. Thus $k\ne0$ forces $g\notin\row D$, and $\rank E=8$. On $k=0$, the factors $4h+3$, $55h-24$, and $3h-5$ are nonzero. Formula~\eqref{eq:cubicres}, already its first entry, shows that the cubic compatibility holds precisely at $h=1/5$. At $G$, $D_2$ has rank three, so the vanishing determinant also gives quadratic compatibility. Hence $\rank E=7$. At $O$, $D$ has rank six but $z_1\ne0$, so the extra row increases its rank to seven. At all other points, $\rank E=8$. Substitution in~\eqref{eq:cubiclocal} proves the dimension formula and the stated local/global distinction.
\end{proof}

\subsection{Proof for $r=2$}\label{app:quinticproof}

The four interior vertices are $v_0,v_1,v_2,Q$. Their direction counts are $(4,4,3,3)$ throughout the domain: the new slope at $v_0$ is $(k+1)/h<-1/2$, different from its old slopes $\infty,\pm1/2$, and at $v_1$ it lies strictly between $-1/2$ and $1$, different from the old slopes $-1,1,-1/2$. The other two counts are immediate. Hence the local ranks are $15,15,14,14$. There are fifteen interior edges, so
\begin{equation}\label{eq:quinticdim}
\dim S_5^2=111-\rank C,\qquad L_2=53.
\end{equation}
Remove the three radial edges and the rows at $Q$, and denote the resulting $45\times72$ matrix by $C_0$.
\subsubsection{The constant matrix}\label{app:quinticmatrix}

Write $U$ in three fifteen-column vertex blocks. Its first three rows have a single nonzero block each, at vertices $v_0,v_1,v_2$, respectively:
\[
(0,-128,0,96,0_{11}),\quad
(192,-32,-64,96,0_{11}),\quad
(-192,-32,64,96,0_{11}).
\]
The fourth row has blocks
\begin{align*}
&(0,-66,0,0,0,0,42,0,-9,0,0,0,-54,0,135),\\
&(132,0,-44,0,96,0,-32,0,96,0,0,0,0,0,0),\\
&(-132,0,44,0,96,0,-32,0,96,0,0,0,0,0,0).
\end{align*}
These four rows are independent and satisfy $UC_0=0$ by direct substitution. Thus $\rank C_0\le41$. The minor with row set
\[
\{0,\ldots,44\}\setminus\{3,18,33,38\}
\]
and column set
\[
\{0,\ldots,17\}\cup\{24,\ldots,35\}\cup\{48,\ldots,58\}
\]
has absolute determinant $2^{125}3^{12}11\ne0$. Therefore $\rank C_0=41$ and these rows give the full left kernel.

\subsubsection{The radial kernel}
The radial block $D$ is $15\times18$. The three lines at $Q$ are distinct, so Lemma~\ref{lem:local} gives $\rank D=3+5+6=14$.
\begin{lemma}\label{lem:radial}
For $i=0,1,8$, set $\ell_i=(y_i-k)(x-h)+(h-x_i)(y-k)$ and
\[
a_0=-2(h-k+2),\quad a_1=3h,\quad a_8=-(h+2k+2),\qquad
u=a_0\ell_0,\ v=a_1\ell_1,\ w=a_8\ell_8.
\]
Then $u+v+w=0$. In the component order $(0,1,8)$, the coefficient vectors of the following four polynomial triples form a basis of $\ker D$:
\begin{gather}
N=(S_1,xS_1,yS_1,S_2),\label{eq:N}\\
S_1=(a_0^3(u+2v),-a_1^3(2u+v),a_8^3(u-v)),\nonumber\\
S_2=(a_0^3(uv+3v^2),a_1^3(3u^2+uv),a_8^3uv).\nonumber
\end{gather}
\end{lemma}
\begin{proof}
The $a_i$ are the signed two-by-two determinants of the other two line-coefficient vectors, proving $u+v+w=0$. They are nonzero on $\mathcal D_2$, and $u,v$ are independent coordinates centered at $Q$. Substituting $w=-u-v$ verifies
\begin{align*}
u^3(u+2v)-v^3(2u+v)+w^3(u-v)&=0,\\
u^3(uv+3v^2)+v^3(3u^2+uv)+w^3uv&=0.
\end{align*}
Thus all four columns lie in the kernel. The first three are independent: an identity $(\lambda_0+\lambda_1x+\lambda_2y)S_1=0$ forces the affine factor to vanish since $S_1$ has a nonzero polynomial component and $\R[x,y]$ is an integral domain. If $S_2$ belonged to their span, then $S_2=fS_1$ for an affine $f$; the third component would imply $uv=f(u-v)$, impossible since $u-v$ does not divide $uv$. There are four independent columns and $\dim\ker D=18-14=4$.
\end{proof}

\subsubsection{A single linear obstruction}
We now reduce the remaining matrix and identify the linear equation that characterizes its rank drop.

By~\eqref{eq:reduction}, $\rank C=55+\rank E$, where $E=UBN$. Put
\begin{gather*}
t=h-k+2,\quad s=h+2k+2,\quad F=13h-5k+10,\\
A_0=10h^2-9hk+20k^2-20k-40,\\
B_0=37h^2+81hk+144h-70k^2+70k+140.
\end{gather*}
\begin{lemma}\label{lem:quinticcalc}
The third row of $E$ is zero. Dividing its first, second, and fourth rows respectively by
\begin{align*}
R_1&=4608h^2(h-2k-2)t^3s^2,\\
R_2&=10368h^3t^2(h+k+2)s^2,\\
R_4&=144h^2t^2s^2
\end{align*}
and omitting the zero row gives
\begin{equation}\label{eq:Ehat}
\widehat E=\begin{pmatrix}
1&0&-1&-\frac92hs\\
1&-2&0&3ts\\
2A_0&-36h(h+2)&B_0&\frac32hts(95h+140k+242)
\end{pmatrix}.
\end{equation}
All these divisions are valid on $\mathcal D_2$.
\end{lemma}
\label{app:quinticcalc}
\begin{proof}[Proof of Lemma~\ref{lem:quinticcalc}]
Let $U_{j,i}$ be row $j$ of $U$ restricted to the block at $v_i$, with $i=0,1,2$. Regard it as a linear functional on the coefficients of degrees $3,4,5$ in centered coordinates $(X,Y)$. Index the rows of $U$ by $j=1,\ldots,4$ here. Let $N_{i\nu}(x,y)$ be the cofactor on edge $(i,9)$ in column $\nu$ of~\eqref{eq:N}, where $i\in\{0,1,8\}$ and $1\le\nu\le4$. Define
\[
n_{i\nu,ab}=[X^aY^b]N_{i\nu}(x_i+X,y_i+Y),\quad
A_i=y_i-k,\quad B_i=h-x_i.
\]
The coefficient in row $(\alpha,\beta)$ of the component at $v_i$ is
\[
\sum_{a,b} n_{i\nu,ab}
\binom{3}{\alpha-a}A_i^{\alpha-a}B_i^{\beta-b},
\qquad \alpha+\beta-a-b=3,
\]
where the sum retains only nonnegative exponents and $0\le\alpha-a\le3$. Contract this expression with $U_{j,i}$ and sum over $i=0,1$:
\[
E_{j\nu}=\sum_{i=0,1}\ \sum_{\substack{\alpha,\beta\ge0\\3\le\alpha+\beta\le5}}
(U_{j,i})_{\alpha\beta}
\sum_{\substack{a,b\ge0\\a+b=\alpha+\beta-3}}
n_{i\nu,ab}\binom{3}{\alpha-a}A_i^{\alpha-a}B_i^{\beta-b}.
\]
The same nonnegative-exponent convention applies to the inner sum. The edge ending at $v_8$ has no old interior-vertex row; the block at $v_2$ receives no radial edge. Thus the third row of $E$ is zero. Substituting the four columns of~\eqref{eq:N} in the finite sum gives, entry by entry,
\[
E_{1,*}=R_1\left(1,0,-1,-\frac92hs\right),\qquad
E_{2,*}=R_2\left(1,-2,0,3ts\right),
\]
\[
E_{4,*}=R_4\left(2A_0,-36h(h+2),B_0,
\frac32hts(95h+140k+242)\right).
\]
These identities yield the stated normalized rows once the row factors are shown to be nonzero.

On $\mathcal D_2$, $h<0$, $t,s>0$, $h+k+2>1+h/2>0$, and $h-2k-2<2h<0$. Each factor $R_1,R_2,R_4$ is therefore nonzero, and division gives~\eqref{eq:Ehat} throughout the domain.
\end{proof}
\begin{proof}[Proof of Theorem~\ref{thm:quintic}]
Denote the rows of~\eqref{eq:Ehat} by $r_1,r_2,r_4$. Set $\beta=18h(h+2)$ and $\alpha=2A_0-\beta$. Direct subtraction gives
\begin{equation}\label{eq:quinticlast}
r_4-\alpha r_1-\beta r_2
=3sF\left(0,0,1,\frac h2(5h+4k+10)\right).
\end{equation}
The first two rows are independent. A linear combination of them with its first two entries zero must itself be zero. Since $s>0$, the displayed residual therefore adds one to the rank precisely when $F\ne0$. Hence $\rank E=2$ if $F=0$, and $3$ otherwise. Formula~\eqref{eq:quinticdim} gives the result. Intersecting $k=(13h+10)/5$ with the strict inequalities defining $\mathcal D_2$ shows that the exceptional set is exactly the open segment joining $(-30/31,-16/31)$ and $(0,2)$.
\end{proof}

\clearpage
\section{Proof of the higher-order dimension formula}\label{app:high}

Each of the $r-1$ insertions adds two faces, three interior edges, and one interior vertex. Thus
\[
f_2=2r+12,\quad f_1^\circ=3r+14,\quad f_0^\circ=r+3.
\]
At $A$, the new upper slope $k/(1+h)$ lies in $(0,1)$, distinct from the four old directions and the negative lower slopes. The corresponding slope at $B$ lies in $(-1,0)$. Hence $n_A=n_B=r+3$. All other $r+1$ interior vertices have three directions, throughout $\mathcal D$. Lemma~\ref{lem:local} therefore gives
\[
\sigma_A=\sigma_B=0,\qquad \sigma_v=\left\lfloor\frac{r^2}{4}\right\rfloor\quad(v\ne A,B).
\]
Substituting in~\eqref{eq:L} yields
\begin{equation}\label{eq:highL}
L_r=11\binom{r+2}{2}+(r+1)\left\lfloor\frac{r^2}{4}\right\rfloor.
\end{equation}

\subsection{The bound $\delta_r\le1$}\label{app:bound}
We bound $\delta_r$ by restricting a left relation to its components at $A$ and $B$. After removing local relations, this restriction is injective; the lower-chain equations then bound its image by a one-dimensional space.

\subsubsection{Dual relations and the image of the fixed-edge matrix}
We express the global defect through dual relations. We then determine the image of the fixed-edge matrix and show that a global relation is determined, modulo local relations, by its two endpoint components.

The defect measures the difference between the sum of the local image dimensions and the rank of $C_\Delta$. We compute this difference by passing to the dual spaces: left relations detect the linear conditions on the image, while local relations account for those already imposed at individual vertices. To write these relations explicitly, identify the dual of $S_d$ with $\mathcal H=\R[U,V,W]_d$ using the differential pairing
\[
\langle f,F\rangle=f(\partial_U,\partial_V,\partial_W)F|_0,
\qquad \langle\ell^m q,F\rangle=\langle q,\ell(\partial)^mF\rangle.
\]
In the monomial bases this pairing is diagonal with nonzero factorial entries, hence nondegenerate. Its second identity turns multiplication by an edge equation into a directional derivative. We use $C_\Delta^{\mathsf T}$ for the dual map under this pairing; in dual bases it is represented by the transpose matrix. Since the two endpoint contributions of an interior edge have opposite signs, a tuple $(F_v)$ belongs to $\ker C_\Delta^{\mathsf T}$ if and only if
\begin{equation}\label{eq:dualcondition}
\ell_e(\partial)^m(F_v-F_{v'})=0
\end{equation}
for every interior edge $e=vv'$, with boundary components set to zero. Such a tuple is called a left relation. Define $J_v^\perp=\{F\in\mathcal H:\langle f,F\rangle=0\ \text{for all }f\in J_v\}$. A tuple supported at $v$ with component in $J_v^\perp$ satisfies all the edge equations by itself; we call it a local relation. Therefore $N_{\loc}=\bigoplus_vJ_v^\perp\subseteq\ker C_\Delta^{\mathsf T}$. Nondegeneracy gives
\[
\begin{aligned}
\dim\ker C_\Delta^{\mathsf T}&=f_0^\circ\dim S_d-\rank C_\Delta,\\
\dim N_{\loc}&=f_0^\circ\dim S_d-\sum_v\dim J_v.
\end{aligned}
\]
Subtracting these identities yields
\begin{equation}\label{eq:dualquotient}
\delta_r=\sum_v\dim J_v-\rank C_\Delta
=\dim\bigl(\ker C_\Delta^{\mathsf T}/N_{\loc}\bigr).
\end{equation}
Thus it suffices to determine the left relations modulo local ones. We shall also use $L_a=a_xU+a_yV+W$ for $a=(a_x,a_y,1)$, for which $\langle q,L_a^r\rangle=r!q(a)$.

Let $E_{\rm old}^{\circ}$ be the interior edges of $\Delta_{\rm old}$. At the same degrees $m=r+1$ and $d=2r+1$, define the old cofactor map by restricting~\eqref{eq:C} to these edges and the four old interior vertices:
\[
C_{\rm old}:\bigoplus_{e\in E_{\rm old}^{\circ}}S_r\longrightarrow
\bigoplus_{v\in\{A,B,X,Y\}}S_d,\qquad
(C_{\rm old}q)_v=\sum_{\substack{e\in E_{\rm old}^{\circ}\\e\ni v}}\varepsilon_{ve}\ell_e^m q_e.
\]
Write $J_v^{\rm old}=\sum_{e\in E_{\rm old}^{\circ},\,e\ni v}\ell_e^mS_r$ for its local image at $v$. In the final triangulation, $C_{\rm old}$ is the block formed by the old-edge columns, with zero rows added for the new vertices. The superscript $\mathrm{old}$ on $J_v^{\rm old}$ thus specifies the incident edges used to form the image.

\begin{lemma}[Independent realization of the fixed-edge vertex components]\label{lem:old}
\[
\im C_{\rm old}=J_A^{\rm old}\oplus J_B^{\rm old}\oplus J_X^{\rm old}\oplus J_Y^{\rm old}.
\]
\end{lemma}
\label{app:old}
We work in affine coordinates: homogenization identifies every polynomial of degree at most $d$ with its degree-$d$ homogeneous form. All ideals below are truncated at degree $d$, and their generators have degree $m$, so their coefficients have degree at most $d-m=r$.

\begin{proof}[Proof of Lemma~\ref{lem:old}]
The inclusion $\im C_{\rm old}\subseteq\bigoplus_vJ_v^{\rm old}$ follows from the definition of each vertex image. For the reverse inclusion, fix arbitrary components $f_v\in J_v^{\rm old}$ for $v=A,B,X,Y$. We construct edge cofactors whose image is exactly $(f_A,f_B,f_X,f_Y)$.

\medskip\noindent\textit{Step 1: components supplied by edges from interior vertices to boundary vertices.}
An interior edge joining an interior vertex to a boundary vertex has one interior endpoint and one boundary endpoint, so its column contributes to only one interior-vertex row block. The following table lists edges representing all the required directions; nonzero constant factors in edge equations can be absorbed into the cofactors.
\[
\begin{array}{c|c|c}
v&\text{boundary endpoints}&\text{edge equations}\\\hline
A&v_1,v_0,v_7&x+1,\ x+1-y,\ x+1+y\\
B&v_2,v_3,v_4&x-1,\ x-1+y,\ x-1-y\\
X&v_1,v_2&x-y-1,\ x+y+1\\
Y&v_5,v_6&x-y+1,\ x+y-1
\end{array}
\]
The two directions at $X$ and the two at $Y$ are already represented by these edges. Thus every component in $J_X^{\rm old}$ or $J_Y^{\rm old}$ can be supplied without affecting the other interior vertices. At $A,B$, put
\[
I_A=((x+1)^m,(x+1-y)^m,(x+1+y)^m),\qquad
I_B=((x-1)^m,(x-1-y)^m,(x-1+y)^m).
\]
The remaining old direction at each endpoint is horizontal. Consequently
\[
J_A^{\rm old}=I_A+(y^m),\qquad J_B^{\rm old}=I_B+(y^m).
\]
In particular, there are $p_A,p_B\in\Pa_r$ such that
\[
f_A\equiv y^mp_A\pmod{I_A},\qquad f_B\equiv y^mp_B\pmod{I_B}.
\]

\medskip\noindent\textit{Step 2: the finite Taylor data that determine these classes.}
Write $u=x+1$, $v=y$ at $A$ and let $I=(u^m,(u-v)^m,(u+v)^m)$ in $\R[u,v]$. Its generators are homogeneous in these centered coordinates. Lemma~\ref{lem:local}, with three directions, gives
\[
\dim(\R[u,v]/I)_{m+j}
=\max\{m+j+1-3(j+1),0\}
=\max\{r-2j-1,0\}.
\]
Set $K_0=\left\lfloor\frac{r-2}{2}\right\rfloor$. For every $K_0<j\le r$, this dimension is zero, so every homogeneous polynomial of degree $m+j$ belongs to $I$. If two degree-at-most-$r$ polynomials have the same Taylor terms of total degree at most $K_0$ at $A$, their difference is a sum of homogeneous terms of degrees $j>K_0$. Multiplication by $y^m=v^m$ therefore puts their difference in $I_A$. The same argument applies at $B$.

More explicitly, write the Taylor expansion in the centered coordinates as $p=\sum_{j=0}^r p_j(u,v)$, where $p_j$ is homogeneous of degree $j$, and set $T_A^{K_0}p=\sum_{j=0}^{K_0}p_j$. Then
\[
v^m(p-T_A^{K_0}p)=\sum_{j=K_0+1}^{r}v^mp_j\in I.
\]
Each summand has degree $m+j$, and the quotient in that degree is zero. Membership can be realized with coefficients of degree $j\le r$ because the generators are homogeneous of degree $m$. Thus discarding the higher Taylor terms is valid within the prescribed cofactor degree bound.

\medskip\noindent\textit{Step 3: one cofactor realizing both endpoint classes.}
Choose the incidence signs on $AB$ to be $+1$ at $A$ and $-1$ at $B$, and use the edge equation $y=0$. Seek its cofactor in the form
\[
q(x,y)=\sum_{j=0}^{K_0}y^jH_j(x).
\]
By Step 2, it is enough to impose
\[
T_A^{K_0}q=T_A^{K_0}p_A,\qquad
T_B^{K_0}q=-T_B^{K_0}p_B.
\]
The minus sign at $B$ is required because the contribution of $AB$ there is $-y^mq$. Group the prescribed Taylor data by powers of $y$:
\[
T_A^{K_0}p_A=\sum_{j=0}^{K_0}y^ja_j(x+1),\qquad
-T_B^{K_0}p_B=\sum_{j=0}^{K_0}y^jb_j(x-1).
\]
Explicitly,
\[
a_j(u)=\sum_{i=0}^{K_0-j}\frac{\partial_x^i\partial_y^jp_A(-1,0)}{i!j!}u^i,
\qquad
b_j(u)=-\sum_{i=0}^{K_0-j}\frac{\partial_x^i\partial_y^jp_B(1,0)}{i!j!}u^i.
\]
Matching the coefficient of each $y^j$ through order $K_0-j$ in $x$ matches exactly all mixed derivatives of total order at most $K_0$. Thus the two-variable problem reduces to independent two-point interpolation problems in $x$.

For each $j$, set $n_j=K_0-j$. There is a unique polynomial $H_j(x)$ of degree at most $2n_j+1$ satisfying
\[
H_j^{(i)}(-1)=a_j^{(i)}(0),\qquad
H_j^{(i)}(1)=b_j^{(i)}(0)\qquad(0\le i\le n_j).
\]
Indeed, the map taking a polynomial of degree at most $2n_j+1$ to its derivatives of orders $0,\ldots,n_j$ at $-1$ and $1$ has equal source and target dimensions $2n_j+2$. Its kernel is zero: a polynomial in the kernel has zeros of multiplicity at least $n_j+1$ at both $-1$ and $1$, and is therefore divisible by $(x+1)^{n_j+1}(x-1)^{n_j+1}$, whose degree is $2n_j+2$. The map is therefore bijective. This proves the required two-point Hermite interpolation directly.
Each term has total degree at most $j+2(K_0-j)+1=2K_0-j+1\le2K_0+1\le r$. Moreover, $H_j(x)-a_j(x+1)$ vanishes to order at least $n_j+1$ at $x=-1$, and $H_j(x)-b_j(x-1)$ does so at $x=1$. Multiplication by $y^j$ raises the total vanishing order at the respective endpoints to at least $j+n_j+1=K_0+1$. Thus $q$ has precisely the prescribed Taylor data and lies in $\Pa_r$. By Step 2,
\[
y^mq\equiv f_A\pmod{I_A},\qquad
-y^mq\equiv f_B\pmod{I_B}.
\]

\medskip\noindent\textit{Step 4: completion by edges from interior vertices to boundary vertices.}
Put this $q$ on $AB$. The residual components $f_A-y^mq\in I_A$ and $f_B+y^mq\in I_B$ can be supplied by the three directions of edges from interior vertices to boundary vertices at their respective endpoints. Such representations use cofactors of degree at most $r$: decompose each residual into homogeneous parts in centered coordinates and express its degree-$n$ part using the degree-$m$ generators with coefficients of degree $n-m\le r$. Finally supply $f_X,f_Y$ by their edges from interior vertices to boundary vertices as in Step 1, and set all unused edge cofactors to zero. No correction changes another interior-vertex component. The resulting image is the prescribed tuple, proving the reverse inclusion.
\end{proof}

The endpoint restriction is the linear map
\begin{equation}\label{eq:endpointmap}
\begin{gathered}
\mathcal R:\ker C_\Delta^{\mathsf T}/N_{\loc}
\longrightarrow (\mathcal H/J_A^\perp)\oplus(\mathcal H/J_B^\perp),
\\
 [(F_v)_v]\longmapsto([F_A],[F_B]).
\end{gathered}
\end{equation}
Here $J_A,J_B$ refer to the final triangulation. Adding a local relation changes $F_A,F_B$ only by elements of $J_A^\perp,J_B^\perp$, respectively, so the map is well defined. Its kernel consists of global relation classes whose endpoint components are local. Injectivity therefore means that if the two endpoint components are local, then every component is local. The following lemma supplies the propagation step needed to prove this assertion; Proposition~\ref{prop:endpoints} also gives canonical representatives for the endpoint classes.
\begin{lemma}\label{lem:twonode}
Fix $a,b>0$ and $m=r+1$, $d=2r+1$. For each nonzero height $t$, let $F_t$ be a homogeneous form of degree $d$ satisfying
\[
(t\partial_U-a\partial_V+at\partial_W)^mF_t=0,\qquad
(t\partial_U+b\partial_V-bt\partial_W)^mF_t=0.
\]
These are the conditions for the two side edges from $(0,t)$ to $(-a,0)$ and $(b,0)$. If $s,t$ are distinct and nonzero and $\partial_U^mF_s=\partial_U^mF_t$, then
\[
\partial_U^mF_s=\partial_U^mF_t=0.
\]
\end{lemma}
\label{app:axis}
\begin{proof}[Proof of Lemma~\ref{lem:twonode}]
Put $L_t=W+tV$. The side equations are $tx-ay+atz$ and $tx+by-btz$, whose operators become $t\partial_U-a\partial_V$ and $t\partial_U+b\partial_V$ in coordinates $(U,V,L_t)$. Their coefficient determinant is $t(a+b)\ne0$. Put $\xi=U+(t/a)V$ and $\eta=U-(t/b)V$. The first derivative is a nonzero multiple of $\partial_\eta$ and the second a nonzero multiple of $\partial_\xi$. Expand $F_t=\sum_\ell f_\ell(U,V)L_t^{d-\ell}$. Each binary layer satisfies $\partial_\xi^mf_\ell=\partial_\eta^mf_\ell=0$, so every monomial $\xi^i\eta^{\ell-i}$ occurring in it has both exponents at most $r$. For $\ell=m+j$, both exponents are consequently at least $\ell-r=j+1$. Thus $f_\ell=0$ for $\ell>2r$, and for $\ell=m+j$, $0\le j\le r-1$,
\[
f_{m+j}=\left[\left(U+\frac taV\right)\left(U-\frac tbV\right)\right]^{j+1}H_j(U,V),
\qquad\deg H_j=r-j-1.
\]
For $\ell<m$, the polynomial $f_\ell$ has degree below $m$ and $\partial_UL_t=0$, so $\partial_U^m(f_\ell L_t^{d-\ell})=0$. Hence the common derivative $G$ has an expansion
\[
G=\sum_{j=0}^{r-1}g_j(U,V)L_t^{r-j},\qquad g_j=\partial_U^mf_{m+j}.
\]
In particular $L_t\mid G$. We claim
\begin{equation}\label{eq:degreeU}
g_j\ne0\quad\Longrightarrow\quad\deg_Ug_j\ge2j-r+1.
\end{equation}
If $2j-r+1\le0$, the bound is automatic. Otherwise let $h=r-j-1<j$ and write $H_j=\sum_{i=0}^hh_iU^{h-i}V^i$. Since $g_j\ne0$, also $H_j\ne0$; let $i_0$ be the smallest index with $h_{i_0}\ne0$. The factor multiplying $H_j$ in $f_{m+j}$ is monic in $U$. Consequently the coefficient of $U^{m+j-i_0}V^{i_0}$ in $f_{m+j}$ is exactly $h_{i_0}$: all terms involving $h_i$ with $i<i_0$ vanish. Differentiating gives
\[
[U^{j-i_0}V^{i_0}]g_j
=\frac{(m+j-i_0)!}{(j-i_0)!}h_{i_0}\ne0.
\]
Since $i_0\le h<j$, the factorials are well defined and $\deg_Ug_j\ge j-i_0\ge j-h=2j-r+1$, proving~\eqref{eq:degreeU}.

The same expansion at $s$ gives $L_s\mid G$. Since $s\ne t$, the linear forms $L_s,L_t$ are relatively prime, so $L_sL_t\mid G$. Suppose $(L_sL_t)^\nu\mid G$ for some $\nu\ge1$. If $2\nu>r$, then $G=0$, since $G$ is homogeneous of degree $r$. Otherwise write $G=(L_sL_t)^\nu K$, where $\deg K=r-2\nu$. Both factors are independent of $U$, hence $\deg_UG\le r-2\nu$.

Expand again in powers of $L_t$. Divisibility by $L_t^\nu$ makes all lower coefficients zero. The coefficient of $L_t^\nu$ is $g_{r-\nu}$, and replacing $W$ by $L_t-tV$ preserves the bound on its $U$-degree. If this coefficient were nonzero,~\eqref{eq:degreeU} would give
\[
r-2\nu\ \ge\ \deg_Ug_{r-\nu}\ \ge\ 2(r-\nu)-r+1=r-2\nu+1,
\]
a contradiction. Thus this coefficient also vanishes, and $L_t^{\nu+1}\mid G$. Interchanging $s,t$ gives $L_s^{\nu+1}\mid G$, hence $(L_sL_t)^{\nu+1}\mid G$. Induction eventually gives a divisor of degree greater than $r$, forcing $G=0$.
\end{proof}

Every old edge is retained and has zero entries at the new vertices. Therefore Lemma~\ref{lem:old} gives $\langle f,F_v\rangle=0$ for every $f\in J_v^{\rm old}$ and $v\in\{A,B,X,Y\}$. Set
\[
D(U,V)=UV(U^2-V^2),\quad L_A=W-U,\quad L_B=W+U,
\quad J=\left\lfloor\frac{r-3}{3}\right\rfloor,\quad h_j=r-3j-3.
\]
\begin{proposition}[Normalized endpoints]\label{prop:endpoints}
Modulo local relations, each left relation has uniquely determined endpoint components of the form
\begin{equation}\label{eq:endpoint}
F_A^{\mathrm{norm}}=\sum_{j=0}^JD^{j+1}H_{A,j}L_A^{r-j},\qquad
F_B^{\mathrm{norm}}=\sum_{j=0}^JD^{j+1}H_{B,j}L_B^{r-j},
\end{equation}
where $H_{A,j},H_{B,j}$ are binary forms of degree $h_j$. The map $\mathcal R$ in~\eqref{eq:endpointmap}, expressed using these representatives, is injective throughout $\mathcal D$.
\end{proposition}
\begin{proof}
Write $F_A=\sum_{\ell=0}^d f_\ell(U,V)L_A^{d-\ell}$, where each $f_\ell$ is binary homogeneous of degree $\ell$. For an edge through $A$, its differential operator kills $L_A$ and therefore acts separately on each binary layer. We first record the resulting divisibility condition.

For any nonzero binary directional derivative $D_e$, choose independent linear coordinates $s,t$ with $D_es=1$ and $D_et=0$. Then $D_e=\partial_s$. If $f$ has degree $m+j$, the equation $D_e^mf=0$ says that every monomial $s^it^{m+j-i}$ has $i\le m-1$. Hence its exponent of $t$ is at least $j+1$, and $t^{j+1}\mid f$. The converse follows by the same monomial calculation. Distinct edge directions give nonproportional factors $t$.

By Lemma~\ref{lem:old}, the component $F_A$ is orthogonal to the four old edge images. Their invariant linear factors are $U,V,U-V,U+V$. Therefore $D^{j+1}\mid f_{m+j}$. The remaining degree is $m+j-4(j+1)=h_j$, so these layers have the form in~\eqref{eq:endpoint}, and vanish for $j>J$.

It remains to identify which layers can be removed as local relations. If $\ell<m$, every $m$th directional derivative of $f_\ell$ is zero. Since all incident-edge operators kill $L_A$, a component $f_\ell L_A^{d-\ell}$ supported only at $A$ satisfies every edge equation and belongs to $J_A^\perp$. Conversely, the final triangulation has $r+3$ distinct directions at $A$. A local layer of degree $m+j$ would have to be divisible by $r+3$ pairwise coprime linear factors, each to the power $j+1$. But
\[
(r+3)(j+1)>m+j\qquad(j\ge0),
\]
so every such layer is zero. Thus $J_A^\perp$ consists exactly of the layers of degree below $m$. Removing them gives the unique representative $F_A^{\mathrm{norm}}$. The argument at $B$ is identical.

For injectivity suppose these endpoints vanish, and subtract local representatives so that $F_A=F_B=0$. The equations on the edges joining a lower vertex to $A$ or $B$ now have zero right-hand sides, as required in Lemma~\ref{lem:twonode}. Along the lower path, the edge equations make all $\partial_U^mF_v$ equal. Lemma~\ref{lem:twonode} applied to any two distinct lower heights forces the common derivative to vanish. Thus $\ell_e(\partial)^mF_v=0$ for every lower vertex $v$ and every incident edge $e$, so each such component is a local relation.

For the upper path put $o=h/(1-k)\in(-1,1)$ and change coordinates to $\xi=x+oy-o$, $\eta=y$. Then $Y,Q$ become $(0,1),(0,k)$ and $A,B$ become $(-(1+o),0),(1-o,0)$. Homogeneously the dual change is $\widetilde F(S,T,R)=F(S,T+oS,R-oS)$. The chain rule gives
\[
\widetilde\ell(\partial)^m\widetilde F
=\bigl(\ell(\partial)^mF\bigr)(S,T+oS,R-oS).
\]
Thus the equations $\ell(\partial)^mF=0$ and the spaces of local relations are preserved. The line $YQ$ is $(1-k)\xi=0$, so its condition equates the two $\partial_S^m$ derivatives. Lemma~\ref{lem:twonode}, with $a=1+o>0$, $b=1-o>0$ and heights $1,k$, makes this derivative zero. Both upper components are local as well. The quotient class is therefore zero.
\end{proof}

From now on, $F_A,F_B$ denote the normalized representatives $F_A^{\mathrm{norm}},F_B^{\mathrm{norm}}$.

\subsubsection{Triangular reparameterization}
We introduce coordinates adapted to an axis vertex and reparameterize $(H_{A,j},H_{B,j})_{j=0}^J$ by an invertible triangular transformation.

Fix a nonzero axis height $t$ and introduce auxiliary coordinates
\[
z=1/t,\quad L=W+tV,\quad P=U+tV,\quad Q'=U-tV.
\]
On $L=0$, one has $L_A=-P$, $L_B=Q'$, $U=(P+Q')/2$, and $V=z(P-Q')/2$. Write the component at $(0,t)$ in the coordinates $(P,Q',L)$ as
\[
F_t=\sum_{\ell=0}^{d}f_{t,\ell}(P,Q')L^{d-\ell},
\qquad f_t:=f_{t,d}=F_t|_{L=0},
\]
where $f_{t,\ell}$ is homogeneous of degree $\ell$ in $P,Q'$. Thus $f_t$ is the binary layer of highest degree, namely the coefficient of $L^0$. The side-edge operators towards $A$ and $B$ become $2t\partial_{Q'}$ and $2t\partial_P$, respectively, and both kill $L$. Applying~\eqref{eq:dualcondition} and comparing the coefficients of $L^0$ gives
\[
\partial_{Q'}^m\bigl(f_t-F_A|_{L=0}\bigr)=0,\qquad
\partial_P^m\bigl(f_t-F_B|_{L=0}\bigr)=0.
\]
The first equation fixes the coefficient of every monomial whose $Q'$ exponent is at least $m$, and the second fixes those whose $P$ exponent is at least $m$. Since $d=2m-1$, every degree-$d$ monomial belongs to exactly one of these two groups. Thus the two equations determine all coefficients of $f_t$, giving
\begin{equation}\label{eq:truncation}
f_t=\Pi_{\deg_{Q'}\ge m}(F_A|_{L=0})+\Pi_{\deg_P\ge m}(F_B|_{L=0}).
\end{equation}
Here $\Pi$ retains the indicated monomials. Inserting~\eqref{eq:endpoint} yields
\begin{equation}\label{eq:canonical}
f_t=\sum_{q=0}^J\bigl[A_q(z)P^{r-q}(Q')^{r+1+q}
+B_q(z)P^{r+1+q}(Q')^{r-q}\bigr].
\end{equation}
\begin{lemma}\label{lem:parameter}
Writing $\ph(z)=z(1-z^2)$, one has
\[
A_q=\ph^{q+1}a_q,\qquad B_q=\ph^{q+1}b_q,\qquad
\deg a_q,\deg b_q\le r-3q-3.
\]
The linear transformation from $(H_{A,j},H_{B,j})_{j=0}^J$ to $(a_q,b_q)_{q=0}^J$ is invertible and triangular.
\end{lemma}
\begin{proof}
Only layers $j\ge q$ contribute to $A_q$. Their contributions are the coefficients of $P^{j-q}(Q')^{m+q}$ in $D^{j+1}H_{A,j}$ after substitution, with the sign from $L_A^{r-j}$. Moreover
\[
D\left(\frac{P+Q'}2,\frac{z(P-Q')}2\right)
=-\frac{\ph(z)}{16}(Q')^4+PE(P,Q',z).
\]
To obtain $P$-degree $j-q$, at least $q+1$ factors $D$ must contribute their $P$-constant terms, proving divisibility by $\ph^{q+1}$. Since $D^{j+1}H_{A,j}$ has degree $m+j$ and contains $U^{j+1}$, its $V$-degree, hence its $z$-degree after substitution, is at most $r$. This proves the quotient degree bound. On the diagonal $j=q$, the contributions are
\[
A_q^{[q]}=(-1)^{r+1}2^{-(m+q)}\ph^{q+1}H_{A,q}(1,-z),\qquad
B_q^{[q]}=2^{-(m+q)}\ph^{q+1}H_{B,q}(1,z).
\]
For a binary form $H(U,V)=\sum_{i=0}^{h_q}c_iU^{h_q-i}V^i$, evaluation at $(1,\pm z)$ gives $\sum_i(\pm1)^ic_iz^i$. Thus each diagonal map is a coefficientwise isomorphism onto all polynomials of degree at most $h_q$. Starting at $q=J$, recover $H_{A,J},H_{B,J}$, subtract their contributions to lower indices, and repeat. This gives an explicit inverse of the triangular transformation.
The divisibility and degree arguments for $B_q$ are the same. 
\end{proof}

\subsubsection{Necessary interpolation conditions}
The edge equations along the lower chain impose interpolation conditions on $(a_q,b_q)_{q=0}^J$. We derive these conditions and the parity of their coefficients.

The $r-1$ lower inverse heights, including $X$, are
\[
\zeta_0=-1,\qquad \zeta_{i+1}=1/t_i\quad(0\le i\le r-3).
\]
Every edge between consecutive vertices of the lower chain lies on $x=0$. Their equations~\eqref{eq:dualcondition} make $\partial_U^mF_v$ the same for all $v\in\{X,P_0,\ldots,P_{r-3}\}$; denote it by $G_-$. In the following formulas, $f_t$ is expressed in $U,V$ by substituting $P=U+tV$, $Q'=U-tV$. Since $L=W+tV$ is independent of $U$, restriction to $L=0$ commutes with $\partial_U^m$, giving
\[
\partial_U^mf_t=G_-(U,V,-tV).
\]
For $0\le\ell\le r$, define
\begin{gather}
K_\ell(q)=[s^\ell](1+s)^{r-q}(1-s)^{r+1+q},\qquad
c_\ell=\frac{(2r+1-\ell)!}{(r-\ell)!},\nonumber\\
T_\ell(z)=\sum_{q=0}^JK_\ell(q)\ph(z)^{q+1}
\bigl(a_q(z)+(-1)^\ell b_q(z)\bigr).\label{eq:T}
\end{gather}
Then $\deg T_\ell\le r$ and expansion of~\eqref{eq:canonical} gives
\[
[U^{r-\ell}V^\ell]\partial_U^mf_t=c_\ell z^{-\ell}T_\ell(z).
\]
If $G_-=\sum_{a+b+c=r}g_{abc}U^aV^bW^c$, the same coefficient on the right is $\sum_{c=0}^{\ell}g_{r-\ell,\ell-c,c}(-t)^c$. Thus
\begin{equation}\label{eq:interpolation}
\begin{gathered}
T_\ell(\zeta_i)=Q_\ell(\zeta_i)\qquad(0\le i\le r-2),\\
Q_\ell(z)=\frac1{c_\ell}\sum_{c=0}^{\ell}(-1)^c g_{r-\ell,\ell-c,c}z^{\ell-c},
\qquad \deg Q_\ell\le\ell.
\end{gathered}
\end{equation}
Only $\ell\le2J+1\le r$ will be used.

\begin{lemma}[Parity of the weights]\label{lem:weights}
Set $x_q=2q+1$, $\mu_q=x_q^2$. There are polynomials $\mathcal E_a,\mathcal O_a$, both of degree exactly $a$ for $0\le a\le J$, such that
\[
K_{2a}(q)=\mathcal E_a(\mu_q),\qquad K_{2a+1}(q)=x_q\mathcal O_a(\mu_q).
\]
\end{lemma}
\begin{proof}
As a formal series the generating function is
\[
(1-s^2)^{(2r+1)/2}\exp\left(\frac{x_q}{2}\log\frac{1-s}{1+s}\right).
\]
Here $\frac12\log\frac{1-s}{1+s}=-s-s^3/3-s^5/5-\cdots$ has zero constant term. In the exponential, the term of degree $j$ in $x_q$ is divisible by $s^j$, so the coefficient of $s^\ell$ has degree at most $\ell$ in $x_q$. Its degree-$\ell$ term comes only from $(-x_qs)^\ell/\ell!$ and the constant term of the prefactor.
Its $s^\ell$ coefficient has degree $\ell$ in $x_q$, with leading coefficient $(-1)^\ell/\ell!$. Replacing $x_q$ by $-x_q$ replaces $s$ by $-s$, proving the parity and degrees.
\end{proof}

\subsubsection{Weighted interpolation and block elimination}
We convert the nodal conditions into a block system. A weighted interpolation determinant gives full column rank for its diagonal blocks, and descending elimination determines its entire kernel.

\begin{lemma}[Weighted interpolation]\label{lem:weighted}
Let $W(x)=\kappa x^a\prod_{j=1}^b(x+\beta_j)$, where $\kappa,\beta_j>0$ and $a,b$ are nonnegative integers. If $p\ge1$ and $a\le\ell\le a+b$ are integers, then the evaluation matrix of
\[
1,x,\ldots,x^{\ell-1},\quad W,xW,\ldots,x^{p-1}W
\]
at any $\ell+p$ strictly increasing positive nodes has positive determinant, with rows in node order and columns in the displayed order. The first group is empty when $\ell=0$.
\end{lemma}
\label{app:weighted}
\begin{proof}[Proof of Lemma~\ref{lem:weighted}]
We first establish positivity of generalized Vandermonde determinants. If $0<u_1<\cdots<u_N$ and $0\le e_1<\cdots<e_N$ are integers, then
\begin{equation}\label{eq:gvpositive}
\det(u_i^{e_j})_{i,j=1}^N>0.
\end{equation}
For $N=1$ this is immediate. Factor $u_i^{e_1}>0$ from row $i$ and write $d_j=e_j-e_1$, so $d_1=0$. Subtract row $i-1$ from row $i$ in descending order $i=N,\ldots,2$, and expand along the first column. The fundamental theorem of calculus and multilinearity give
\[
\frac{\det(u_i^{e_j})}{\prod_{i=1}^Nu_i^{e_1}}
=\left(\prod_{j=2}^Nd_j\right)
\int_{u_1}^{u_2}\!\cdots\!\int_{u_{N-1}}^{u_N}
\det(s_i^{d_j-1})_{\substack{1\le i<N\\2\le j\le N}}
\,ds_{N-1}\cdots ds_1.
\]
The exponents $d_2-1<\cdots<d_N-1$ are nonnegative integers. On the interior of the integration box, $0<s_1<\cdots<s_{N-1}$, so the integrand is positive by induction. By continuity it is nonnegative on the boundary. All intervals have positive length and all $d_j>0$, proving~\eqref{eq:gvpositive}.

Now put $D=a+b$, $N=\ell+p$, and write $W=\sum_kw_kx^k$, taking $w_k=0$ outside $0\le k\le D$. Form the coefficient matrix $A$ of the polynomials in Lemma~\ref{lem:weighted}, with rows numbered $0,\ldots,D+p-1$. Its first $\ell$ columns are the unit vectors $e_0,\ldots,e_{\ell-1}$, and its last $p$ columns have entries $w_{i-j}$ for $0\le j<p$. If $V$ is the matrix of the powers $0,\ldots,D+p-1$ evaluated at the $N$ nodes, the required evaluation matrix is $VA$.

Let $S$ be the lower shift on indices $0,\ldots,D+p-1$, so $S_{i,j}=1$ for $i=j+1$ and is zero otherwise. The square convolution matrix is
\[
T=(w_{i-j})=\kappa S^a\prod_{j=1}^b(\beta_jI+S).
\]
Every minor of $S$ and of $\beta I+S$ is nonnegative. Indeed, in a determinant term a row index can only equal its matched column index or exceed it by one. An inversion would match $i<i'$ with $j>j'$, which requires $i\ge j\ge j'+1$ and $i'\le j'+1$, contradicting $i<i'$. Thus every nonzero term uses the order-preserving permutation, has positive sign, and is a product of nonnegative entries. Cauchy--Binet now shows that every minor of $T$ is nonnegative. The finite truncation is exact because these factors only increase indices, so an index cannot leave the range and later return.

An $N$-row maximal minor of $A$ is zero unless it contains all rows $0,\ldots,\ell-1$. If it does, expansion along the unit columns gives a minor of $T$ with positive sign. Hence every maximal minor of $A$ is nonnegative. For the row set $0,\ldots,N-1$, the remaining minor of $T$ has rows $\ell,\ldots,\ell+p-1$ and columns $0,\ldots,p-1$. It is strictly positive: in the product for $T$, let $S^a$ shift all selected indices by $a$, choose exactly $\ell-a$ of the $b$ bidiagonal factors to shift all indices by one, and use the diagonal entries in every other factor. This is a valid strictly positive term in the repeated Cauchy--Binet expansion because $0\le\ell-a\le b$.

Finally, every maximal minor of $V$ is positive by~\eqref{eq:gvpositive}. In
\[
\det(VA)=\sum_{\substack{I\subseteq\{0,\ldots,D+p-1\}\\|I|=N}}
\det V_{:,I}\det A_{I,:},
\]
all terms are nonnegative and the term with $I=\{0,\ldots,N-1\}$ is positive. This proves the lemma, including $\ell=0$.
\end{proof}

For $0\le q\le J$, put $n=r-2$, $p_q=r-3q-2=n-3q\ge1$, and define
\[
\xi_i=-1-\frac1{t_i}\quad(0\le i<n),\qquad
W_q(x)=x^q(x+1)^{q+1}(x+2)^{q+1}.
\]
The ordering of the heights gives $0<\xi_0<\cdots<\xi_{n-1}$. For $p\ge1$ let $V_p=(\xi_i^j)_{0\le i<n,\,0\le j<p}$, and let $V_0$ have no columns. Set
\[
E_q=\operatorname{diag}(W_q(\xi_0),\ldots,W_q(\xi_{n-1}))V_{p_q}.
\]
Thus, for a coefficient vector $v=(v_0,\ldots,v_{p_q-1})^{\mathsf T}$ and $P(x)=\sum_jv_jx^j$,
\[
E_qv=\begin{pmatrix}W_q(\xi_0)P(\xi_0)\\\vdots\\W_q(\xi_{n-1})P(\xi_{n-1})\end{pmatrix}.
\]
\par

\begin{figure}[htbp]\centering
\begin{tikzpicture}[x=1.35cm,y=1cm,font=\small,>=stealth]
\draw[->](-2.6,0)--(4.6,0)node[right]{$x$};
\foreach \x in {-2,-1,0}{\draw[accent,thick](\x,-.09)--(\x,.09);\node[below]at(\x,-.1){$\x$};}
\foreach \x/\lab in {.65/\xi_0,1.4/\xi_1,3.9/\xi_{n-1}}{\fill[warm](\x,0)circle(2pt);\node[below]at(\x,-.1){$\lab$};}
\node at(2.65,0){$\cdots$};
\node[accent,above]at(-1,.35){$W_q=x^q(x+1)^{q+1}(x+2)^{q+1}$};
\node[warm,above]at(2.6,.35){$0<\xi_0<\cdots<\xi_{n-1}$};
\end{tikzpicture}
\caption{The nonpositive roots of the weight and the positive evaluation nodes. The root at $0$ has multiplicity $q$ (and is absent for $q=0$); the roots at $-1,-2$ have multiplicity $q+1$.}\label{fig:ordered}
\end{figure}
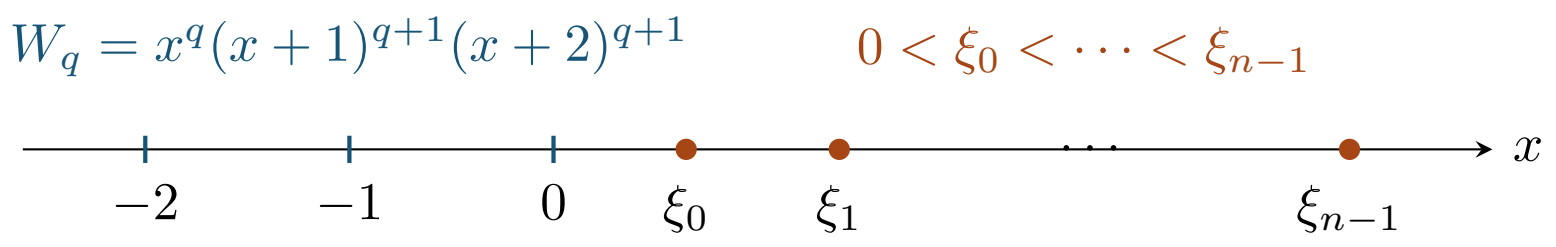

\par
For $0<\ell<n$, define $R_\ell:\R^n\to\R^{n-\ell}$ by subtracting the degree-less-than-$\ell$ interpolant on the first $\ell$ nodes from the remaining data:
\[
(R_\ell v)_{i-\ell}=v_i-\sum_{j=0}^{\ell-1}
\left(\prod_{\substack{0\le s<\ell\\s\ne j}}
\frac{\xi_i-\xi_s}{\xi_j-\xi_s}\right)v_j,
\qquad \ell\le i<n.
\]
Set $R_0=I_n$; for $\ell\ge n$, let $R_\ell$ have no rows. For every $\ell\ge0$, $\ker R_\ell=\im V_\ell$, where $\im V_\ell=\R^n$ for $\ell\ge n$. In particular, these kernels are nested.

Let $\mathbf c_q,\mathbf e_q$ be the coefficient vectors, in increasing powers of $x$, of $a_q(-1-x)+b_q(-1-x)$ and $a_q(-1-x)-b_q(-1-x)$, respectively. Since
\[
\ph(-1-x)=x(x+1)(x+2),\qquad
\frac{\ph(-1-x)^{q+1}}x=W_q(x),
\]
the interpolation conditions~\eqref{eq:interpolation} become
\begin{equation}\label{eq:orderedblocks}
\begin{aligned}
\sum_{q=0}^JK_{2a}(q)R_{2a}E_q\mathbf c_q&=0,\\
\sum_{q=0}^JK_{2a+1}(q)R_{2a+1}E_q\mathbf e_q&=0
\qquad(0\le a\le J).
\end{aligned}
\end{equation}
Indeed, at $\zeta_0=-1$ all $T_\ell$ vanish, so $Q_\ell(-1)=0$. For $\ell\ge1$, the polynomial $Q_\ell(-1-x)/x$ has degree less than $\ell$. At the remaining nodes, division by $\xi_i>0$ gives precisely the corresponding equation with $R_\ell$. Conversely, such a low-degree interpolant, multiplied by $x$ and substituted back, recovers $Q_\ell$ including its value at $-1$. For $\ell=0$, $Q_0=0$ and the same statement uses $R_0=I_n$.

Write $\mathbf1=(1,\ldots,1)^{\mathsf T}\in\R^n$.
\begin{proposition}\label{prop:one}
The full solution space of~\eqref{eq:orderedblocks} is
\[
\mathbf c_q=0\ (0\le q\le J),\qquad
\mathbf e_q=0\ (1\le q\le J),\qquad
\mathbf e_0=\rho E_0^{-1}\mathbf1\quad(\rho\in\R).
\]
Consequently the endpoint pairs satisfying these lower-chain necessary conditions form a one-dimensional space.
\end{proposition}
\begin{proof}
First, if $p_a+\ell\le n$ and $a\le\ell\le3a+2$, then $R_\ell E_a$ has full column rank. To see this, suppose $R_\ell E_av=0$. There is a polynomial $P$ of degree less than $\ell$ whose values equal $E_av$. On any $p_a+\ell$ nodes, Lemma~\ref{lem:weighted}, applied with $b=2a+2$, makes $[V_\ell,E_a]$ invertible. Thus both the coefficients of $P$ and $v$ vanish. In particular,
\[
\rank(R_{2a}E_a)=p_a\quad(0\le a\le J),\qquad
\rank(R_{2a+1}E_a)=p_a\quad(1\le a\le J),
\]
because $p_a+2a=n-a$ and $p_a+2a+1=n-a+1$.

Next put $\mu_q=(2q+1)^2$ and $\pi_a(\mu)=\prod_{j=0}^{a-1}(\mu-\mu_j)$, with $\pi_0=1$. By Lemma~\ref{lem:weights}, expand $\pi_a$ in the bases $\mathcal E_0,\ldots,\mathcal E_a$ and $\mathcal O_0,\ldots,\mathcal O_a$. For each $j\le a$, the even equation at level $j$ says that $\sum_qK_{2j}(q)E_q\mathbf c_q$ belongs to $\im V_{2j}\subseteq\im V_{2a}$. Apply $R_{2a}$ to this vector and then take the indicated constant linear combination. The odd equations are treated with $\im V_{2j+1}\subseteq\im V_{2a+1}$. Since $\pi_a(\mu_q)=0$ for $q<a$, this gives
\[
\sum_{q=a}^J\pi_a(\mu_q)R_{2a}E_q\mathbf c_q=0,
\qquad
\sum_{q=a}^J(2q+1)\pi_a(\mu_q)R_{2a+1}E_q\mathbf e_q=0.
\]
At $a=J$, the even equation has only its diagonal term, whose matrix has full column rank and whose scalar $\pi_J(\mu_J)$ is nonzero; hence $\mathbf c_J=0$. Repeating for $a=J-1,\ldots,0$ gives all $\mathbf c_q=0$. The odd equations give $\mathbf e_J=\cdots=\mathbf e_1=0$ by the same argument. The remaining equation is $R_1E_0\mathbf e_0=0$. Here
\[
E_0=\operatorname{diag}((\xi_i+1)(\xi_i+2))V_n
\]
is invertible. Since $\ker R_1=\Span\{\mathbf1\}$, this yields the stated formula for $\mathbf e_0$. Conversely, these vectors satisfy every original equation: the even sums vanish, and each odd sum is a scalar multiple of $\mathbf1$, which belongs to every $\ker R_{2a+1}$. This also covers $r=3$, where $n=1$ and $R_1$ has no rows. Finally, Lemma~\ref{lem:parameter} and the invertible changes $z=-1-x$ and $(a_q,b_q)\mapsto(\mathbf c_q,\mathbf e_q)$ transfer this one-dimensional kernel to the endpoint pairs.
\end{proof}
Together with Proposition~\ref{prop:endpoints}, this proves $\delta_r(h,k)\le1$ throughout $\mathcal D$ for every ordered lower chain.

\subsection{The bound $\delta_r\ge1$ when $h=0$}\label{sec:obstruction}
For $h=0$, we construct a left relation with a nonzero class modulo local relations. The construction starts from an interpolation polynomial $p$, which determines the endpoint components and their extension to all vertices.
\subsubsection{A polynomial with positive coefficients}
We construct the interpolation polynomial $p$ used in $F_A^*,F_B^*$ below and prove positivity of its coefficients.

\begin{lemma}\label{lem:positive}
The unique polynomial $p$ of degree at most $r-3$ satisfying
\begin{equation}\label{eq:p}
(t_i-1)p(t_i)=-t_i^{r-1}\qquad(0\le i\le r-3)
\end{equation}
has strictly positive coefficients.
\end{lemma}
\begin{proof}
Put $n=r-2$ and $K(t)=\prod_{i=0}^{r-3}(t-t_i)=\sum_{j=0}^ne_jt^{n-j}$, where $e_0=1$ and all $e_j>0$. With $S=K(1)$ and $c_0=1-S^{-1}$, define
\[
p(t)=\frac{(t-c_0)K(t)-t^{n+1}}{t-1}.
\]
The numerator vanishes at $1$ and has degree at most $n$, so this is a polynomial of degree at most $n-1$. It satisfies~\eqref{eq:p}, proving existence; uniqueness follows from the $n$ distinct nodes. The coefficient of $t^{n-\ell}$ is
\[
b_\ell=e_\ell-S^{-1}\sum_{j=\ell}^ne_j\qquad(1\le\ell\le n).
\]
To obtain this coefficient formula, compare the coefficient of $t^{n+1-\ell}$ in $(t-1)p(t)=(t-c_0)K(t)-t^{n+1}$. It gives $b_\ell-b_{\ell-1}=e_\ell-c_0e_{\ell-1}$, with $b_0=0$. Summing and using $S=\sum_{j=0}^ne_j$ yields the displayed expression.
The elementary symmetric functions of the positive numbers $-t_i$ satisfy $e_\ell e_j\ge e_{\ell+j}$ for $\ell+j\le n$: the product contains every square-free term of the latter with positive multiplicity. The remaining products are strictly positive. Hence
\[
Sb_\ell=e_\ell\sum_{j=1}^ne_j-\sum_{j=\ell+1}^ne_j
>\sum_{j=1}^{n-\ell}e_{\ell+j}-\sum_{j=\ell+1}^ne_j=0.
\]
The strict inequality also covers $\ell=n$: the second sum is then empty, whereas $e_n\sum_{j=1}^ne_j>0$. For $\ell<n$, the terms with $j>n-\ell$ on the left are positive and have no counterpart on the right. This includes $n=1$, corresponding to $r=3$.
\end{proof}
Define the homogeneous polynomial $H(U,V)=(-V)^{r-3}p(-U/V)$ by expansion, including at $V=0$, and set
\begin{equation}\label{eq:star}
F_A^*=DH(U,V)L_A^r,\qquad
F_B^*=(-1)^rDH(U,-V)L_B^r.
\end{equation}

\begin{proposition}[Axis extension]\label{prop:axisextension}
For $h=0$ and every $0<k<1$, the endpoint pair~\eqref{eq:star} extends to a nonlocal left relation. In particular $\delta_r(0,k)=1$, and the fixed lower-chain endpoint space is spanned by this pair.
\end{proposition}
\begin{proof}
Take boundary components zero and put
\begin{align}
F_X=F_{P_i}&=2^{-r}(W-V)(W^2-U^2)^r,\nonumber\\
F_Y=F_Q&=-2^{-r}\frac{(k+1)p(k)}{k^{r-1}}(W+V)(W^2-U^2)^r.\label{eq:axisextension}
\end{align}

\begin{figure}[htbp]
\centering
\begin{tikzpicture}[x=2.3cm,y=2cm,font=\small]
\coordinate(A)at(-1,0);\coordinate(B)at(1,0);
\coordinate(X)at(0,-1);\coordinate(Y)at(0,1);
\coordinate(Q)at(0,.5);\coordinate(P0)at(0,-.75);\coordinate(P1)at(0,-.3);
\draw[gray] (A)--(B);
\foreach \p in {Y,Q}{\draw[warm] (A)--(\p)--(B);}
\foreach \p in {X,P0,P1}{\draw[accent] (A)--(\p)--(B);}
\draw[warm,very thick] (Y)--(Q);
\draw[accent,very thick] (X)--(P0)--(P1);
\foreach \p in {A,B,X,Y,Q,P0,P1}{\fill (\p)circle(1.5pt);}
\node[left,align=center]at(A){$A$\\$F_A^*$};
\node[right,align=center]at(B){$B$\\$F_B^*$};
\node[above]at(Y){$Y$};\node[right,fill=white,inner sep=1pt]at(Q){$Q=(0,k)$};
\node[below]at(X){$X$};
\node[right,fill=white,inner sep=1pt]at(P0){$P_0$};
\node[right,fill=white,inner sep=1pt]at(P1){$P_{r-3}$};
\node[fill=white,inner sep=0pt]at(0,-.52){$\vdots$};
\draw[warm] (0,.78)--(1.45,.78) node[right]{$F_Y=F_Q=F_+$};
\draw[accent] (0,-.92)--(1.45,-.92) node[right]{$F_X=F_{P_i}=F_-$};
\end{tikzpicture}
\[
F_-=2^{-r}(W-V)(W^2-U^2)^r,\qquad
F_+=C_k(W+V)(W^2-U^2)^r,
\]
\[
C_k=-2^{-r}\frac{(k+1)p(k)}{k^{r-1}}.
\]
\caption{Extension of the endpoint pair on the axis. All lower-chain vertices carry $F_-$, and $Y,Q$ carry $F_+$. Thus the component difference on each thick axis edge is zero. The side edges satisfy the derivative identities in the proof. Boundary components are zero. The lower chain is schematic; for $r=3$ it has just one vertex $P_0$.}\label{fig:extension}
\end{figure}

Figure~\ref{fig:extension} shows the vertex components of this extension.
To check all edge conditions, let $a_t=tx-y+tz$, $b_t=tx+y-tz$, and $\gamma(t)=-t(t^2-1)p(t)$. Since $a_t(\partial)L_A=0$ and $b_t(\partial)L_B=0$,
\[
a_t(\partial)^mF_A^*=m!\gamma(t)L_A^r,\qquad
b_t(\partial)^mF_B^*=(-1)^m m!\gamma(t)L_B^r.
\]
For $F=(uW+vV)(W^2-U^2)^r$, the product rule gives
\[
a_t(\partial)^mF=m!2^rt^r(ut-v)L_A^r,\qquad
b_t(\partial)^mF=(-1)^m m!2^rt^r(ut-v)L_B^r.
\]
On the lower path $u=2^{-r}$, $v=-2^{-r}$. Equation~\eqref{eq:p} matches the endpoint values at $t=t_i$; at $t=-1$ both sides vanish. On the upper path $u=v=C_k=-2^{-r}(k+1)p(k)/k^{r-1}$, and $2^rC_kk^r(k-1)=\gamma(k)$; at $t=1$ both sides vanish. This checks the side edges, including their collinear edges from interior vertices to boundary vertices.

The axis edges have identical components at their two endpoints, so their differences vanish. At $A,B$, the four old directional derivatives vanish because $D$ contains the four required linear factors. These checks cover every interior edge of the mesh. The relation is nonlocal since
\[
a_{t_i}(\partial)^mF_A^*=m!t_i^r(t_i+1)L_A^r\ne0,
\]
so the component at $A$ is not a local relation. Thus $\delta_r(0,k)\ge1$, and the opposite inequality follows from Appendix~\ref{app:bound}. The exhibited nonzero endpoint pair also attains the one-dimensional bound in Proposition~\ref{prop:one}.
\end{proof}
The conditions used in Proposition~\ref{prop:one} involve only the fixed lower chain and are independent of $(h,k)$. The nonzero pair~\eqref{eq:star} satisfies these conditions, so it spans their one-dimensional solution space. Thus, for every $(h,k)\in\mathcal D$, the normalized endpoints of any left relation have the form
\begin{equation}\label{eq:scalar}
(F_A,F_B)=\lambda(F_A^*,F_B^*).
\end{equation}

\subsection{Vanishing of $\delta_r$ when $h\ne0$}\label{app:offaxis}
The preceding construction determines every possible endpoint pair up to a scalar. We pair the representative $(F_A^*,F_B^*)$ with a suitable image of the cofactor map and show that the pairing is nonzero when $h\ne0$. Thus a left relation must have zero endpoint class, and injectivity gives $\delta_r=0$.
\subsubsection{An image supported at $A,B$}
We construct a cofactor vector whose image is supported at the two endpoints. Pairing this image with the endpoint components yields the scalar condition for extension to a moving vertex.

Put $c=1-k$ and use the following homogeneous equations for $AQ,BQ,YQ$:
\[
a=k(x+z)-(1+h)y,\quad b=k(x-z)+(1-h)y,\quad L=cx+h(y-z).
\]
With $\alpha=(c-h)/(2k)$ and $\beta=(c+h)/(2k)$, one has $\alpha,\beta>0$ and $L=\alpha a+\beta b$. Define degree-$r$ cofactors
\begin{align}
q_L&=L^r,\nonumber\\*
q_A&=-\sum_{j=m}^{d}\binom dj\alpha^j\beta^{d-j}a^{j-m}b^{d-j},\nonumber\\*
q_B&=-\sum_{j=0}^{r}\binom dj\alpha^j\beta^{d-j}a^jb^{d-j-m}.\label{eq:q}
\end{align}
The exponents are nonnegative. The binomial expansion gives
\[
a^mq_A+b^mq_B+L^d=0.
\]
Orient these edges toward $Q$. The component at $Q$ vanishes; the components at $A,B,Y$ are $-a^mq_A,-b^mq_B,-L^d$. The two edges from interior vertices to boundary vertices from $Y$ to $(1,2),(-1,2)$ have equations $g_A=x-y+z$, $g_B=x+y-z$, and
\[
L=\frac{c-h}{2}g_A+\frac{c+h}{2}g_B.
\]
Splitting this binomial expansion at exponent $m$ expresses $L^d=g_A^mu_A+g_B^mu_B$, with $u_A,u_B\in S_r$. These edges touch no other interior vertex, so they cancel the component at $Y$ without changing the endpoint components. The resulting cofactor vector $\mathcal U(h,k)$ satisfies
\begin{equation}\label{eq:components}
C_{\Delta_r(h,k)}\mathcal U(h,k)=(-a^mq_A,-b^mq_B,0,\ldots,0).
\end{equation}

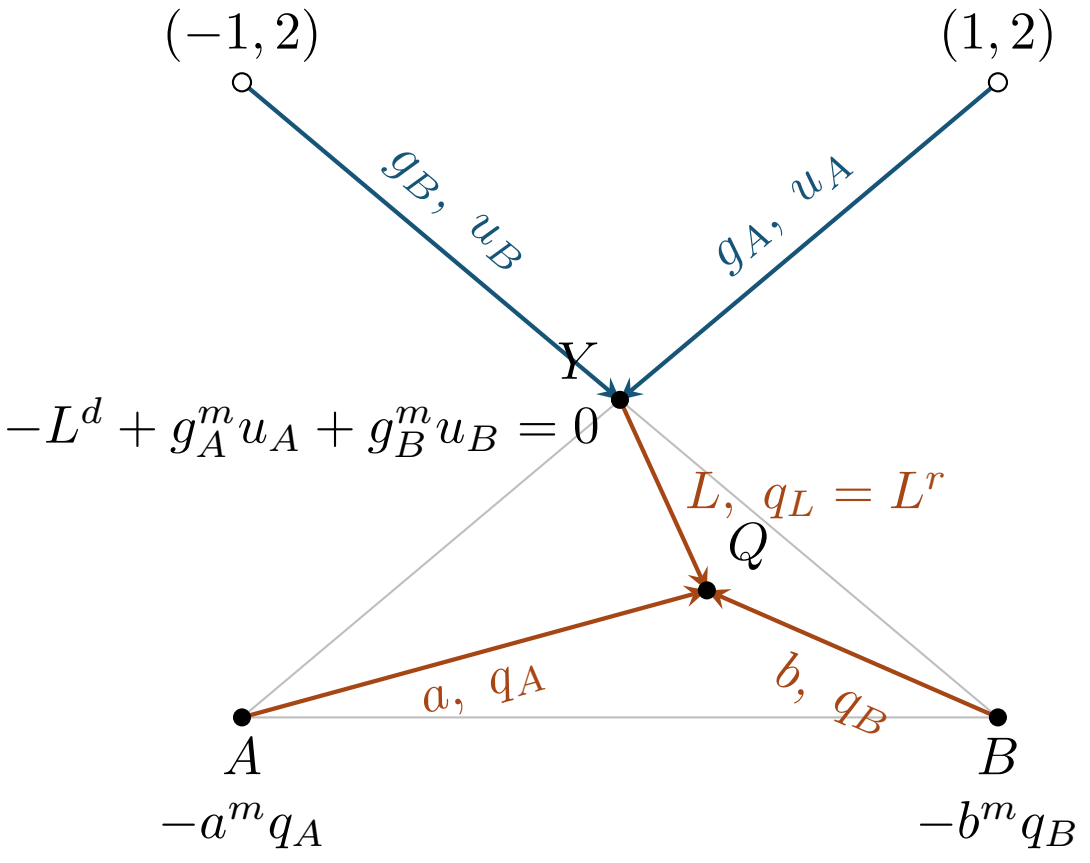
\begin{figure}[htbp]
\centering
\begin{tikzpicture}[x=2.5cm,y=2.1cm,font=\small,>=stealth]
\coordinate(A)at(-1,0);\coordinate(B)at(1,0);
\coordinate(Y)at(0,1);\coordinate(Q)at(.23,.4);
\coordinate(TA)at(1,2);\coordinate(TB)at(-1,2);
\draw[gray!50] (A)--(Y)--(B)--cycle;
\draw[warm,thick,->] (A)--(Q) node[midway,below,sloped]{$a,\ q_A$};
\draw[warm,thick,->] (B)--(Q) node[midway,below,sloped]{$b,\ q_B$};
\draw[warm,thick,->] (Y)--(Q) node[midway,right]{$L,\ q_L=L^r$};
\draw[accent,thick,->] (TA)--(Y) node[midway,above,sloped]{$g_A,\ u_A$};
\draw[accent,thick,->] (TB)--(Y) node[midway,above,sloped]{$g_B,\ u_B$};
\foreach \p in {A,B,Y,Q}{\fill (\p)circle(1.7pt);}
\foreach \p in {TA,TB}{\draw[fill=white] (\p)circle(1.7pt);}
\node[above]at(TA){$(1,2)$};\node[above]at(TB){$(-1,2)$};
\node[below,align=center]at(A){$A$\\$-a^mq_A$};
\node[below,align=center]at(B){$B$\\$-b^mq_B$};
\node[left,align=right]at(Y){$Y$\\$-L^d+g_A^mu_A+g_B^mu_B=0$};
\node[above right]at(Q){$Q$};
\end{tikzpicture}
\[
\text{Component at }Q:\quad
a^mq_A+b^mq_B+L^mq_L=0.
\]
\caption{Cofactors and vertex components. Each edge label gives its linear equation and its cofactor; arrows specify the orientation. The three orange edges cancel at $Q$. The two blue interior edges joining interior vertices to boundary vertices add $L^d$ at $Y$, cancelling its initial component $-L^d$. Only the displayed components at $A,B$ remain in the interior-vertex rows of $C_{\Delta_r(h,k)}\mathcal U$.}\label{fig:loads}
\end{figure}

The orientations, cofactors, and remaining components are shown in Figure~\ref{fig:loads}.\par
Pairing any left relation with this image and using~\eqref{eq:scalar} gives the necessary condition
\begin{equation}\label{eq:necessary}
0=\lambda\Omega_r(h,k),\qquad
\Omega_r=\langle-a^mq_A,F_A^*\rangle+\langle-b^mq_B,F_B^*\rangle.
\end{equation}

\begin{lemma}[Exact obstruction]\label{lem:obstruction}
On $\mathcal D$,
\begin{equation}\label{eq:omega}
\Omega_r(h,k)=\Phi_r(h,k)\,k(c^2-h^2)
\bigl[\mathcal F_k(1+h)-\mathcal F_k(1-h)\bigr],
\end{equation}
where
\[
\Phi_r=(-1)^rm!r!\binom dm(2k)^{r-1}\alpha^r\beta^r\ne0,\qquad
\mathcal F_k(s)=s^{r-2}(s+k)p(k/s).
\]
In particular $\Omega_r(h,k)=0$ if and only if $h=0$.
\end{lemma}
\begin{proof}
Because $a(\partial)L_A=0$ and $DH$ has degree $m$,
\[
\langle a^mq,F_A^*\rangle=m!r!(DH)(k,-(1+h))q(-1,0,1).
\]
Indeed, the product rule places all $m$ derivatives on $DH$, because every term differentiating $L_A^r$ contains $a(\partial)L_A=0$. A homogeneous degree-$m$ polynomial evaluated by its $m$th directional derivative gives $m!$ times its value at that direction. The remaining pairing is $\langle q,L_A^r\rangle=r!q(-1,0,1)$. The same reasoning applies at $B$, retaining the factor $(-1)^r$ in $F_B^*$.
Put $\Gamma_k(s)=ks^{r-2}(s^2-k^2)p(k/s)$. The factors $(DH)(k,-(1+h))$ and $(-1)^rD(k,1-h)H(k,-(1-h))$ evaluate to $\Gamma_k(1+h)$ and $(-1)^{r+1}\Gamma_k(1-h)$. At $A$, $a(A)=0$ and $b(A)=-2k$, so only $j=m$ survives in $q_A$. At $B$, $a(B)=2k$ and $b(B)=0$, so only $j=r$ survives in $q_B$. Thus~\eqref{eq:q} gives:
\[
q_A(A)=-\binom dm\alpha^m\beta^r(-2k)^r,\qquad
q_B(B)=-\binom dm\alpha^r\beta^m(2k)^r.
\]
Keeping the signs of the components in~\eqref{eq:necessary} gives
\begin{align*}
\Omega_r
&=(-1)^rm!r!\binom dm(2k)^r\alpha^r\beta^r
\bigl[\alpha\Gamma_k(1+h)-\beta\Gamma_k(1-h)\bigr]\\
&=\Phi_r\bigl[(c-h)\Gamma_k(1+h)-(c+h)\Gamma_k(1-h)\bigr].
\end{align*}
Since $\Gamma_k(s)=k(s-k)\mathcal F_k(s)$, and $1\pm h-k=c\pm h$, both terms acquire the common factor $k(c-h)(c+h)$. This gives~\eqref{eq:omega}. Write $p(t)=\sum_{j=0}^{r-3}p_jt^j$, where $p_j>0$ by Lemma~\ref{lem:positive}. Then
\[
\mathcal F_k'(s)=\sum_{j=0}^{r-3}p_jk^j
\bigl[(r-1-j)s^{r-2-j}+k(r-2-j)s^{r-3-j}\bigr]>0\quad(s>0).
\]
On $\mathcal D$, $k>0$, $c^2-h^2>0$, and $1\pm h>0$. Strict monotonicity proves that the bracket vanishes precisely at $h=0$.
\end{proof}

\begin{proof}[Proof of Theorem~\ref{thm:high}]
For $h\ne0$, Lemma~\ref{lem:obstruction} and~\eqref{eq:necessary} force $\lambda=0$. Injectivity of $\mathcal R$ in Proposition~\ref{prop:endpoints} then makes every left relation local, so $\delta_r=0$. For $h=0$, Proposition~\ref{prop:axisextension} gives $\delta_r=1$. Equations~\eqref{eq:split} and~\eqref{eq:highL} finish the proof.
\end{proof}

\clearpage



\begin{thebibliography}{99}
\small
\bibitem{Sch84}
L. L. Schumaker, Bounds on the dimension of spaces of multivariate piecewise polynomials,
\emph{Rocky Mountain Journal of Mathematics} 14 (1984), 251--264.

\bibitem{AS87}
P. Alfeld and L. L. Schumaker, The dimension of bivariate spline spaces of smoothness $r$ for degree $d\ge4r+1$,
\emph{Constructive Approximation} 3 (1987), 189--197.

\bibitem{Hong}
D. Hong, Spaces of bivariate spline functions over triangulation,
\emph{Approximation Theory and its Applications} 7 (1991), 56--75.

\bibitem{AS90}
P. Alfeld and L. L. Schumaker, On the dimension of bivariate spline spaces of smoothness $r$ and degree $d=3r+1$,
\emph{Numerische Mathematik} 57 (1990), 651--661.

\bibitem{Universal}
X. Wu and J. Deng, Universal dimension formula for planar $C^r$ splines at degree $3r+1$,
preprint, 2026.

\bibitem{MS75}
J. Morgan and L. R. Scott, The dimension of the space of $C^1$ piecewise polynomials,
unpublished manuscript, 1975.

\bibitem{Diener}
D. Diener, Instability in the dimension of spaces of bivariate piecewise polynomials of degree $2r$ and smoothness order $r$,
\emph{SIAM Journal on Numerical Analysis} 27 (1990), 543--551.

\bibitem{WZD}
X. Wu, Y. Zheng, and J. Deng, A geometric criterion for the dimension of Morgan--Scott spline spaces $S_{2r}^r(\Delta_{\mathrm{MS}})$,
preprint, 2026.

\bibitem{Stiller}
H. K. Schenck and P. F. Stiller, Cohomology vanishing and a problem in approximation theory,
\emph{Manuscripta Mathematica} 107 (2002), 43--58.

\bibitem{YS}
B. Yuan and M. E. Stillman, A counter-example to the Schenck--Stiller \textquotedblleft$2r+1$\textquotedblright\ conjecture,
\emph{Advances in Applied Mathematics} 110 (2019), 33--41.

\bibitem{SSY}
H. Schenck, M. Stillman, and B. Yuan, A new bound for smooth spline spaces,
\emph{Journal of Combinatorial Algebra} 4 (2020), 359--367.

\bibitem{Billera}
L. J. Billera, Homology of smooth splines: generic triangulations and a conjecture of Strang,
\emph{Transactions of the American Mathematical Society} 310 (1988), 325--340.

\bibitem{WD}
X. Wu and J. Deng, Geometry-dependent rank defect in $C^1$ cubic spline space,
preprint, 2026.

\bibitem{Landscape}
X. Wu and J. Deng, The dimension landscape of bivariate $C^2$ splines,
preprint, 2026.

\bibitem{LW}
Z. Luo and R. Wang, Structure and dimension of multivariate spline space of lower degree on arbitrary triangulation,
\emph{Journal of Computational and Applied Mathematics} 195 (2006), 113--133.

\bibitem{SS97}
H. Schenck and M. Stillman, Local cohomology of bivariate splines,
\emph{Journal of Pure and Applied Algebra} 117--118 (1997), 535--548.

\bibitem{Survey}
H. Schenck, Algebraic methods in approximation theory,
\emph{Computer Aided Geometric Design} 45 (2016), 14--31.
\end{thebibliography}
\end{document}